\documentclass[a4paper]{article}

\usepackage{amsmath,amsthm,amssymb,mathtools}
\usepackage{comment} 
\usepackage[T1]{fontenc} 
\usepackage[mathcal]{euscript}
\usepackage{graphicx}
\usepackage[hidelinks]{hyperref}
\usepackage[top=3cm,bottom=3cm,left=2.5cm,right=2.5cm]{geometry}
\usepackage{tikz}
\usepackage[export]{adjustbox}
\usepackage{appendix}
\usepackage{setspace}
\usepackage{tcolorbox}
\numberwithin{equation}{section}
\newtheorem{theorem}{Theorem}[section]

\newtheorem{lemma}[theorem]{Lemma}

\makeatletter
\g@addto@macro\bfseries{\boldmath}
\makeatother

\usetikzlibrary{shapes.geometric,decorations.markings,intersections,calc,positioning}

\tikzset{->-/.style={decoration={
  markings,
  mark=at position .6 with {\arrow{>}}},postaction={decorate}}}

\tikzset{-<-/.style={decoration={
  markings,
  mark=at position .6 with {\arrow{<}}},postaction={decorate}}}

\tikzset{%
    add/.style args={#1 and #2}{
        to path={%
 ($(\tikztostart)!-#1!(\tikztotarget)$)--($(\tikztotarget)!-#2!(\tikztostart)$)%
  \tikztonodes},add/.default={.2 and .2}}
}  

\tikzset{
    extended line/.style={shorten >=-#1,shorten <=-#1},
    extended line/.default=1cm]
}

\tikzset{line through/.style args={#1 parallel to line through #2 and #3 and
length #4}{insert path={%
let \p1=($(#3)-(#2)$),\n1={atan2(\y1,\x1)} in (#1) -- ++ (\n1:#4)}}}

\definecolor{col1}{rgb}{0.4, 0.69, 0.2}%cyan(process)
\definecolor{col2}{rgb}{0.96, 0.29, 0.54}%fulvous

\definecolor{green(ryb)}{rgb}{0.4, 0.69, 0.2}
\definecolor{frenchrose}{rgb}{0.96, 0.29, 0.54}
\definecolor{persianblue}{rgb}{0.11, 0.22, 0.73}
\definecolor{jade}{rgb}{0.0, 0.66, 0.42}
\definecolor{limegreen}{rgb}{0.2, 0.8, 0.2}

\begin{document}
\title{ \vspace{-2.5cm}
Geometric approach for the identification \\ of Hamiltonian systems of quasi-Painlev\'e type}
\author{Marta Dell'Atti{$\,^{1,2}$} \\[-.8ex] \small{\url{m.dell-atti@uw.edu.pl}} \\[1ex] Thomas Kecker{$\,^{1}$} \\[-.8ex]  \small{\url{thomas.kecker@port.ac.uk}} \\[.5ex]
{$\,^{1}$}{ \footnotesize School of Mathematics and Physics, University of Portsmouth, UK}
\\[-.2ex]
{$\,^{2}$}{ \footnotesize Corresponding author}
}
\date{}
\maketitle
\vspace{-.5cm}
\begin{abstract}
   \hspace*{-1ex} Some new Hamiltonian systems of quasi-Painlev\'e type are presented and the  analogue of Okamoto's space of initial conditions computed. Using the geometric approach that was introduced originally for the identification problem of Painlev\'e equations, comparing the irreducible components of the inaccessible divisors arising in the blow-up process, we find bi-rational coordinate changes between some of these systems that give rise to the same global Hamiltonian structure. This scheme thus gives a method for identifying Hamiltonian systems up to bi-rational maps, which is performed in this article for systems of quasi-Painlev\'e type having singularities that are either square-root type algebraic poles or ordinary poles.
\end{abstract}

\vspace*{-1ex} 

\tableofcontents

\newpage
\section{Introduction}

We report on research that is motivated originally by the following question in the theory of complex differential equations: what kinds of singularities can a local analytic solution of a non-linear ordinary differential equation develop when continued in the complex plane? In the case of linear differential equations the situation is clear: apart from a number of fixed singularities (points at which one or several coefficients of the equation become ill-defined), the solutions can be analytically continued along any path in the plane. If fact, the solutions form a linear space, spanned by a set of fundamental solutions.

For non-linear equations, the situation is not at all obvious. Firstly, besides the fixed singularities, a solution can have (usually an infinite number of) movable singularities. Roughly speaking, these are singular points of a solution the locations of which depend on the initial data for the equation: when going from one solution with initial data $(y(z_0),y'(z_0)) = (y_0,\eta_0)$ to a solution with modified initial data $(y_0 + \epsilon,\eta_0 + \zeta)$ with {$|\epsilon|,|\zeta| \ll 1$}, the position of the singularities will change in a continuous fashion. Secondly, for a given non-linear differential equation, it is not at all clear how to describe the space of all solutions. 

While there are abundant examples of non-linear equations with interesting behaviour, in this article we restrict ourselves to a particular case important in mathematical physics, the Hamiltonian systems. These systems of equations are derived from a scalar function $H(x,y;z)$, which would be an integral of motion if not explicitly $z$-dependent. While we consider only Hamiltonians $H$ polynomial in two dependent variables $x(z)$ and $y(z)$, we do allow for an explicit (analytic) $z$-dependence of the coefficients,
\begin{equation}
\label{hamiltonian}
H = H(x(z),y(z);z) = \sum_{(j,k) \in I} a_{jk}(z)\, x(z)^j\, y(z)^k,
\end{equation}

\vspace*{-2ex}

\noindent
where $I \subset \mathbb{N}^2$ is some finite index set. The Hamiltonian equations of motion are given by 
\begin{equation}
\label{ham_eqns}
\begin{cases}
 y' \equiv \dfrac{d y}{d z} = \dfrac{\partial H}{\partial x} = \displaystyle \sum_{(j,k) \in I} j\, a_{jk}(z)\, x(z)^{j-1}\, y(z)^k\,, \\[3ex] x' \equiv \dfrac{d x}{d z} = -\dfrac{\partial H}{\partial y} = - \displaystyle \sum_{(j,k) \in I} k\, a_{jk}(z)\, x(z)^j\, y(z)^{k-1}.
 \end{cases} 
\end{equation}
Thus, starting from a local, analytic solution $x(z),y(z)$ for $z \in U$, $U \subset \mathbb{C}$ a domain, we are interested in what types of singularities the solution can develop when continued under the flow of the Hamiltonian vector field $\left( {\partial H}/{\partial x} , -{\partial H}/{\partial y} \right)$. {Note that, for a non-autonomous system, the flow of the vector field should be considered on a three-dimensional manifold with coordinates $(x,y,z)$, but we interpret the picture on the two-dimensional manifold with coordinates $(x,y)$, with the flow itself being $z$-dependent.} The nature of the singularities is often determined by the highest power terms, but possibly other terms, in the Hamiltonian. E.g., if the coefficient function of $a_{j_1 k_1}(z_\ast)$ of a leading term in the Hamiltonian happens to be zero at a singular point $z_\ast \in \mathbb{C}$ of the solution $(x(z),y(z))$, the nature of that singularity will be different to the generic case when ($a_{j_1 k_1}(z_\ast) \neq 0$). Those points would be classed as fixed singularities. Since in this article we are predominantly studying the movable singularities, we will usually assume the coefficients of the leading-order terms to be non-zero constants. In particular, by redefining $x(z),y(z)$, in a neighbourhood $U$ of a point $z_\ast$ where $a_{j k}(z_\ast) \neq 0$ for some $(j,k) \in \{(j_1,k_1),(j_2,k_2)\}$ one can set these two functions to be non-zero constants, with the remaining coefficients being analytic functions in $U$.

In general, there are many different types of movable singularities that can occur in the solutions of non-linear differential equations, such as poles, algebraic poles, logarithmic and essential singularities, and even non-isolated singularities (e.g.\ encountered in certain equations of degree $\geq 3$, such as the Chazy equation \cite{Chazy1911}). In this paper, the equations we are studying are of such nature that the only movable singularities that can occur, by continuation of the solution along finite-length paths in the $z$-plane, are at most algebraic poles, i.e.\ at a singularity $z_\ast$ the solution can be represented by a series expansion in a fractional power of $z-z_\ast$ (Puiseux series). This property is now known as the quasi-Painlev\'e property, a notion introduced by S. Shimomura in the articles~\!\cite{Shimomura2006,Shimomura2008}. { In particular, the property guarantees a certain regularity of the solutions, in the way that locally, the solutions extend over a Riemann surface with only a finite number of sheets, and singularities with a more complicated branching, such as logarithmic ones, are excluded (globally, however, i.e.\ over the whole complex plane, the picture is much more complicated in general).} It is a generalisation of the notion of equations with the Painlev\'e property, for which all movable singularities of all solutions have to be poles. We briefly discuss this class of equations in the next section but will quickly move on to equations of quasi-Painlev\'e type in section \ref{sec:intro_quasi-Painleve}, which this article is mainly about.%\footnote{Throughout the paper we use the abbreviation `qsi-P' for quasi-Painlev\'e, not to be confused with the $q$-Painlev\'e in the context of discrete Painlev\'e and orthogonal polynomials usually found in the literature. }.

\subsection{Equations of Painlev\'e type}

There exists a class of non-linear second-order differential equations whose solutions have a particularly simple singularity structure. Motivated by a function theoretic viewpoint to find new, transcendental meromorphic functions defined as solutions of differential equations, P. Painlev\'e and his school classified non-linear equations of the form $y''=R\big(y(z),y'(z);z\big)$, $R$ rational in its first two arguments, with the property that all movable singularities of all their solutions are poles~\!\cite{Painleve1900}. Among this class, they isolated roughly $50$ types of equations, of which a family of $6$ equations gives rise to solutions which could not be expressed in terms of formerly known special functions. These equations, known as the six Painlev\'e equations, are
\begin{align}
\label{Painleve1}
    \text{P}_{\text{I}}\colon &\, y'' =  6\,y^2 + z \,,\\[2ex]
\label{Painleve2}
    \text{P}_{\text{I\!I}}\colon &\, y'' =  2\,y^3 + zy + \alpha \,, \\[2ex]
\label{Painleve3}
    \text{P}_{\text{I\!I\!I}}\colon &\, y'' =  \frac{(y')^2}{y} - \frac{y'}{z} + \frac{1}{z}(\alpha\, y^2 + \beta) + \gamma\, y^3 + \frac{\delta}{y}\,, \\[.7ex]
\label{Painleve4}
    \text{P}_{\text{I\!V}}\colon &\, y'' =  \frac{(y')^2}{2\,y} + \frac{3}{2} \,y^3 + 4\,z\,y^2 +2(z^2-\alpha)y + \frac{\beta}{y}\,, \\[.7ex]
\label{Painleve5}
    \text{P}_{\text{V}}\colon &\, y'' =  \frac{3\,y-1}{2\,y(y-1)} (y')^2 - \frac{y'}{z} + \frac{(y-1)^2}{z^2} \left( \alpha y + \frac{\beta}{y} \right) + \frac{\gamma\, y}{z} + \frac{\delta\, y(y+1)}{y-1} \,,\\[.7ex]
\label{Painleve6}
\begin{split} 
    \text{P}_{\text{V\!I}}\colon &\, y'' =  \frac{1}{2} \left( \frac{1}{y} + \frac{1}{y-1} + \frac{1}{y-z} \right) (y')^2 - \left( \frac{1}{z} + \frac{1}{z-1} + \frac{1}{y-z} \right) y' \\[.7ex] & \qquad + \frac{y(y-1)(y-z)}{z^2(z-1)^2} \left( \alpha  + \beta\, \frac{z}{y^2} + \gamma\, \frac{z-1}{(y-1)^2} + \delta\, \frac{z(z-1)}{(y-z)^2} \right)\,,
\end{split} 
\end{align}
where $\alpha,\beta,\gamma,\delta$ are complex parameters. For the equations~\!\eqref{Painleve1},~\!\eqref{Painleve2} and~\!\eqref{Painleve4}, which do not have any fixed singularities, the Painlev\'e property means that every solution can be continued to a globally meromorphic function in $\mathbb{C}$. For the other Painlev\'e equations, which also have fixed singularities, their solutions can be continued to meromorphic functions on a covering surface of the plane punctured at the set of fixed singular points ($\mathbb{C} \setminus \{0\}$ for equations~\!\eqref{Painleve3} and~\!\eqref{Painleve5}, $\mathbb{C} \setminus \{0\,,1\}$ for equation~\!\eqref{Painleve6}). The Painlev\'e equations, although the first ones are strikingly simple, exhibit an exceedingly rich mathematical content in terms of symmetry, underlying algebraic and geometric structure, that cannot be anticipated directly from the form of these equations. 

Each of the six Painlev\'e equations can be written {in Hamiltonian form}, with a Hamiltonian that is polynomial in its dependent variables. These were studied extensively by Okamoto in a series of four papers~\!\cite{okamoto1,okamoto2,okamoto3,okamoto4}, who also studied their B\"acklund transformations: these are discrete symmetries of the equations while acting on the parameter space $(\alpha,\beta,...)$ as affine Weyl groups. Furthermore, for each of the Painlev\'e equations, Okamoto~\!\cite{Okamoto1979} introduced the \textit{space of initial conditions}, assigning a geometric object to each of these equations. This is obtained by first extending the phase space to a compact, rational surface (such as $\mathbb{CP}^2$ or $\mathbb{CP}^1 \times \mathbb{CP}^1$), which is repeatedly \textit{blown up} at the base points, points of indeterminacy of the vector field defining the system in any coordinate chart covering this space. In the blow-up process, certain exceptional curves arise, forming an inaccessible divisor (in the sense of algebraic geometry) of the space. By studying the configuration of the irreducible components of these divisors, Sakai~\!\cite{Sakai2001} classified all differential and discrete Painlev\'e equations according to their inaccessible divisor classes, assigning to each equation a surface type and a symmetry type by which they can be identified. 

The space of initial conditions is essentially the space obtained from the initially compact surface under various cascades of blow-ups, with the inaccessible components of the exceptional divisors removed (which results in a non-compact surface). Using the blow-up procedure outlined above it is possible to construct a global atlas for the space of initial conditions with a canonical $2$-form {for the regularised system free of points of indeterminacy}. Here, in each coordinate chart of the space, it can be achieved that the Hamiltonian is polynomial with respect to the standard symplectic form,

\vspace*{-2ex}

\begin{equation}
\omega = dy \wedge dx = du_\text{f} \wedge dv_\text{f}\,,
\end{equation}
where $(x,y)$ are the original variables of the system and $(u_\text{f},v_\text{f})$ the coordinates after the final blow-up in any cascade of blow-ups. In this way, Takano et al.~\!\cite{Takano97,Takano99} proved that for each Painlev\'e system there exists a unique global Hamiltonian structure.

\subsection{Equations with the quasi-Painlev\'e property}
\label{sec:intro_quasi-Painleve}

It is S. Shimomura~\!\cite{Shimomura2006,Shimomura2008} who first studied classes of equations with what he calls the quasi-Painlev\'e property, such as 
\begin{equation}
\begin{aligned}
    y''(z) &= \frac{2(2k+1)}{(2k-1)^2} y^{2k} + z, \quad k \in \mathbb{N}, \\[1ex]
    y''(z) &= \frac{k+1}{k^2} y^{2k+1} + z y + \alpha, \quad \alpha \in \mathbb{C}, \quad k \in \mathbb{N} \setminus \{2\}.
    \end{aligned}
\end{equation}

This property for an equation or system of equations, being a generalisation of the ordinary Painlev\'e property, demands that all singularities of the equation that can be obtained by continuation of the solution along finite length paths, are at worst algebraic poles (instead of just ordinary poles in the case of the Painlev\'e equations). That is, near a movable singularity $z_\ast \in \mathbb{C}$, the solution is represented by a Puiseux series in a cut neighbourhood of that point,
{
\begin{equation}
	\label{formal_series} 
y(z) = \sum_{j=j_0}^{\infty} c_j (z-z_\ast)^{j/n}, \qquad j_0 \in \mathbb{Z}, n \in \mathbb{N}.
\end{equation}
}
Filipuk and Halburd~\!\cite{halburd1} study more general second-order differential equations with the quasi-Painlev\'e property (though not referred to by this name). There, it was shown for the class of equations
\begin{equation}
\label{Halburd_class}
y''(z) = P\big(y(z);z\big),
\end{equation}
where $P$ is a polynomial in $y$ (of degree $\deg_y P = N$, say) with $z$-analytic coefficients, that under certain {\it resonance conditions} the equation~\eqref{Halburd_class} has the quasi-Painlev\'e property. These conditions arise by inserting a formal series of the form~\eqref{formal_series} as solution for the equation. One can easily obtain the different types of leading-order behaviour for such a solution, $y \sim c_{j_0} (z-z_\ast)^{j_0/n}$. The next step is to set up a recurrence relation for the coefficients $c_j$, $j\geq 1$, of the form $(j-2N-2)c_j = Q_j(c_0,c_1,\dots,c_{j-1})$. When recursively solving for the coefficients $c_j$, there is an obstruction at a certain index $j$ when the numerical factor multiplying $c_j$ vanishes and we cannot solve for $c_j$ for this particular index $j=2N+2$. This is called a resonance, and if the right-hand side of the recursive relation is not identically $0$ (the resonance condition), no formal algebraic series solutions of the form~\!\eqref{formal_series} exists. The resonance conditions, for each possible leading-order behaviour of the solution being satisfied, are thus equivalent with the existence of sufficiently many formal Puiseux series solutions. This is essentially saying that the equation passes the quasi-Painlev\'e test, which is necessary for it to have the quasi-Painlev\'e property. The tricky part in the papers~\!\cite{halburd1,halburd2} (and other papers such as~\!\cite{Shimomura2003} for proofs of the Painlev\'e property for the Painlev\'e equations) is to show that the existence of these formal series solutions is indeed sufficient for an equation to have the quasi-Painlev\'e property, that is to show that there cannot exist any other types of movable singularities in solutions of these equations. For example, for $\deg_y P = 2$ and $\deg_y P = 3$ in~\!\eqref{Halburd_class} we recover the first and second Painlev\'e equations \eqref{Painleve1} and~\!\eqref{Painleve2} respectively, whereas for $\deg_y P = 4$ and $\deg_y P = 5$ we have the following equations,
\begin{align}
\label{quasiP4}	y'' & = y^4 + a_2(z)\, y^2 + a_1(z)\, y + a_0(z)\,, &\qquad &\text{ with }a_2''(z) = 0\,, \\[2ex]
\label{quasiP5} y'' & = y^5 + b_3(z)\, y^3 + b_2(z)\, y^2 + b_1(z)\, y + b_0(z)\,, &\qquad &\text{ with }b_3''(z) = 0,~~ (b_3^2 + 4\,b_1)' = 0\,.
\end{align}
The solutions of these equations have leading-order $y(z) \sim (z-z_\ast)^{-2/3}$ for equation~\!\eqref{quasiP4}, respectively $y(z) \sim \pm (z-z_\ast)^{-1/2}$ for equation~\!\eqref{quasiP5} at any movable singularity $z_\ast \in \mathbb{C}$. Note that, although there is some restriction on the coefficients of these equations, there are arbitrary (analytic) functions, here the coefficients $a_1(z)$, $a_0(z)$ and $b_2(z)$, $b_0(z)$, which are not restricted by the quasi-Painlev\'e property of these equations. One of the authors~\!(TK) generalises these ideas by introducing a class of Hamiltonian systems $H=H\left(q(z),p(z);z\right)$, polynomial in the variables $q,p$, for which the Hamiltonian equations have the quasi-Painlev\'e property~\!\cite{Kecker2016}. The equations~\!\eqref{quasiP4} and~\!\eqref{quasiP5} are special cases of such Hamiltonian systems. 

Filipuk and Kecker~\!\cite{KeckerFilipuk} show that the concept of the Okamoto's space of initial conditions can be extended to the class of quasi-Painlev\'e equations and Hamiltonian systems of quasi-Painlev\'e type, and in~\!\cite{KeckerFilipuk} these spaces are constructed for equations~\!\eqref{quasiP4} and~\!\eqref{quasiP5}, among other Hamiltonian systems. In particular, it is shown that, starting from an equation with arbitrary (analytic) functions as coefficients, through the process of regularising the equations via cascades of blow-ups, one can obtain the conditions on the coefficients under which these equations have the quasi-Painlev\'e property. These conditions are found when imposing the regularity of the system of equations obtained after the final blow-up in each cascade. In general, these conditions can be seen to lead to an additional cancellation taking place in the equations, after which they can be integrated without the need to introduce logarithms in the solution. In contrast to the Painlev\'e case, one has to make an additional change of dependent and independent variables for the system to become regular, as would be expected. This process is essentially equivalent to performing the quasi-Painlev\'e test. Furthermore, by introducing an auxiliary function $W$, one can show that the exceptional divisors introduced by the the blow-ups are inaccessible by the solutions, thus establishing an alternative way of proving the quasi-Painlev\'e property for these systems. The auxiliary function is obtained by adding certain correction terms to the Hamiltonian with the effect that $W$ remains bounded at the movable singularities.

\subsection{Overview of this article}
In this article we consider examples of different types of geometries that can arise as {the analogue of} spaces of initial conditions for certain equations and systems with the quasi-Painlev\'e property. For the purpose of this paper, we restrict ourselves to Hamiltonian systems $H\!\left(x(z),y(z);z\right)$, polynomial in $x,y$, with coefficients analytic in $z$ on some common domain $U \subset \mathbb{C}$. Further, we restrict to the case where the leading coefficients of the Hamiltonian are non-zero, and therefore can be absorbed, by re-scaling, into the definition of the dependent variables $x$, $y$, rendering the leading coefficients constants. A point $z_0 \in \mathbb{C}$ where a leading coefficient vanishes is referred to as a \textit{fixed singularity} of the system, while in this article we are predominantly interested in the behaviour of the solutions at \textit{movable singularities}. 
We mainly consider systems in this article for which the leading-order behaviour at a singularity $z_\ast$ of a solution $\left(x(z),y(z)\right)$ is such that $x(z)$, $y(z)$ are either simple poles or algebraic poles of the form 
\begin{equation}
\label{leading_order}
x(z), y(z) \sim c_{-1} (z-z_\ast)^{-1/2}, 
\end{equation}
but we also encounter some systems which exhibit mixed types of singularities. This is interesting as so far in the literature most examples of equations with the quasi-Painlev\'e property only have one type of movable singularity.

In section~\!\ref{geom_approach}, we outline the geometric approach for the identification of Hamiltonian systems via their Okamoto's spaces of initial conditions. After introducing the basic notions, we apply the approach to the example of the Painlev\'e-II equation in section~\!\ref{P2review}, for which we present various Hamiltonian forms. { While in the Painlev\'e case, the geometric approach allows us to explicitly compute a symplectic, bi-rational change of variables for two systems with isomorphic spaces of initial conditions, in this paper we explore to what extent this method can be extended to the quasi-Painlev\'e case. }

In section~\!\ref{sec:quasi_Painleve}, we consider various Hamiltonian systems with the quasi-Painlev\'e property, which we identify via their surface types. In particular, we consider several Hamiltonian systems giving rise to equation~\!\eqref{quasiP5} and others that can be seen as a generalisation of the Painlev\'e-IV equation~\!\eqref{Painleve4}. While in the Painlev\'e case, { all relevant components of the inaccessible divisor are curves of self-intersection $-2$,} giving rise to Dynkin diagrams under the intersection form, in the cases of quasi-Painlev\'e we consider here, the surface types, i.e.\ the configurations of the irreducible components of the inaccessible divisor of a system after all blow-ups, are such that there is exactly one $-3$-curve, while all other components are $-2$-curves. 
While the equation~\eqref{quasiP5} falls into this class (which we denote of quasi-Painlev\'e-II type), in section~\ref{quasi-PainleveII} we present three different Hamiltonian systems that give rise to this equation (by elimination of $x$) and find explicit relationships between these Hamiltonian systems via their spaces of initial conditions.
Further, in section~\ref{quasi-PainleveIV}, we present Hamiltonian systems that give rise to an equation that we consider to be a quasi-Painlev\'e analogue of the fourth Painlev\'e equation~\!\eqref{Painleve4}. In these systems, $x$ and $y$ can have algebraic poles of the form~\!\eqref{leading_order} or ordinary poles (in which case $x$ has a pole and $y$ a zero). This is { an} example of a system in which mixed types of movable singularities occur.

Using the blow-up procedure, in the case of the Painlev\'e equations one can obtain a symplectic atlas of the space of initial conditions for which the Hamiltonian is polynomial in all final charts $(u_\text{f},v_\text{f})$, whereas the $2$-form is of the canonical form $\omega = dy \wedge dx = d u_\text{f} \wedge d v_\text{f}$. To achieve this, one needs to make certain intermediate changes of variables (other than blow-ups) to ensure that the coefficient of $\omega$ remains explicitly $z$-independent. In the quasi-Painlev\'e case, instead of obtaining a symplectic atlas on which the system is defined (i.e.\ where the $2$-form defining the symplectic structure is $\omega = du_\text{f} \wedge dv_\text{f}$ in any coordinate chart), it was shown in~\!\cite{FilipukStokes} that the $2$-form has a zero in the charts after final blow-ups, 
\begin{equation}
\label{quasi_symplectic_form}
\omega = u_\text{f}^{k-1} \, du_\text{f} \wedge dv_\text{f}, \quad k \in \mathbb{N},
\end{equation}
while the Hamiltonian remains polynomial in all such charts. Again, to achieve this various intermediate changes of variables are required. The integer $k$ then indicates the leading-order behaviour at a movable singularity $z_\ast \in \mathbb{C}$ of a solution in this chart: $u_\text{f}(z) \sim (z-z_\ast)^{-1/k}$. For the different cascades of blow-ups of a given system, the value of $k$ may differ in the various charts after final blow-ups. We will encounter such cases in the following, meaning that the solutions have mixed types of movable singularities. While in the Painlev\'e case the transformations leading from the original system of equations to the one after a final blow-up in a cascade are indeed symplectic (with $k=1$ in equation~\!\eqref{quasi_symplectic_form}), in the quasi-Painlev\'e case the form~\!\eqref{quasi_symplectic_form} of $\omega$, with some $k>1$ is the best one can do.

\section{Constructing spaces of initial conditions}
\label{geom_approach}
In this section we review the process of constructing the Okamoto's space of initial conditions and carry out the procedure in the example of the Painlev\'e-II equation~\!\eqref{Painleve2} in subsection~\!\ref{P2review}. We do this for four different Hamiltonians giving rise to this equation. By comparing the spaces of initial conditions of any two cases, one can find a bi-rational change of variables between the respective Hamiltonian systems. This is what we refer to as the geometric approach to identifying certain Hamiltonian systems. We will then apply this identification procedure in later sections of this article to equations with the quasi-Painlev\'e property. 

The procedure of constructing the space of initial condition for a given system starts by extending the phase space of the dependent variables to a compact rational surface, covered by a finite atlas of coordinate charts, in which we re-write the Hamiltonian equations. Note that the bi-rational changes of variables between the different charts of this atlas, as well as the subsequent changes of variables under blow-ups, are in general non-symplectic, i.e.\ do not preserve the $2$-form $\omega = dy \wedge dx$. An initial choice has to be made regarding the compact surface from which to start the blow-up process. In the literature, frequently $\mathbb{CP}^1 \times \mathbb{CP}^1$ is found when performing blow-ups for the Painlev\'e equations, however, in this article we start from $\mathbb{CP}^2$, as we find that it simplifies the calculations for the comparison of the spaces of initial conditions for two related systems. The difference between the two approaches is that in the $\mathbb{CP}^1 \times \mathbb{CP}^1$ case, the coordinates $x$, $y$ are directly assigned to one of the factor spaces, while in the $\mathbb{CP}^2$ case, all lines are equivalent divisors. In other words, the Picard lattice $\text{Pic}(\mathbb{CP}^1 \times \mathbb{CP}^1) = \mathbb{Z} \times \mathbb{Z}$, whereas for $\mathbb{CP}^2$ we have $\text{Pic}(\mathbb{CP}^2) = \mathbb{Z}$, see also section~\!\ref{picard}. The results obtained through the blow-up process, in particular the conditions on the equations to have the quasi-Painlev\'e property and the nature of the movable singularities, do not depend on the choice of the initial compact surface.

\subsection{Blowing up and blowing down}
We consider a differential system in some coordinate chart, which in general is expressed in terms of rational functions in the variables $(x,y)$,
\begin{equation} \label{eq:dyn_sys}
\begin{aligned}
x' = \frac{P_{1}\big(x(z),y(z);z\big)}{Q_{1}\big(x(z),y(z);z\big)}\,, \qquad y' = \frac{P_{2}\big(x(z),y(z);z\big)}{Q_{2}\big(x(z),y(z);z\big)}\,.
\end{aligned}
\end{equation}
At certain points $(x,y)$, the vector field (right-hand side of the equations) can become indeterminate. For a fixed $z_0 \in \mathbb{C}$, these are the points $(x,y)$ at which $P_{j}\big(x(z_0),y(z_0);z_0\big) = 0 = Q_{j}\big(x(z_0),y(z_0);z_0\big)$ for either $j=1$, $j=2$, or both. It must be noted, however, that not at all points of indeterminacy the solution necessarily breaks down. For example, in the analysis below we sometimes find points at which the vector field is indeterminate in some chart, while in a different chart covering that point it is well-defined and therefore unproblematic. 

\begin{figure}[t]
    \centering
    \includegraphics[width=.35\textwidth]{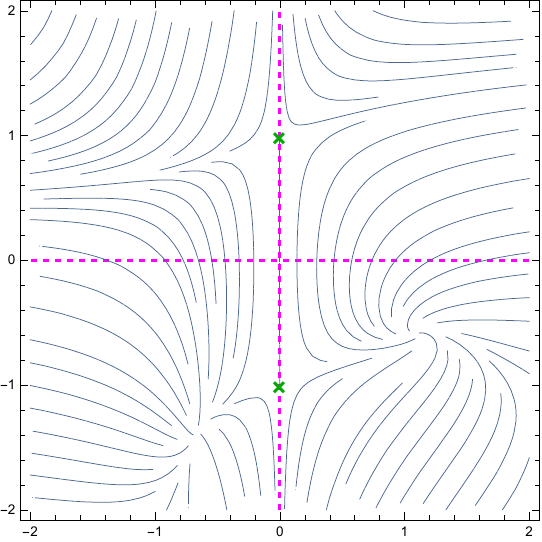} \hspace{8ex} \includegraphics[width=.35\textwidth]{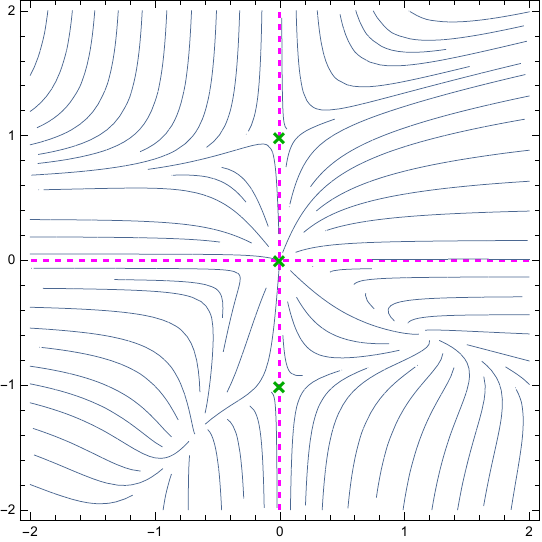} 
    \caption{Geometry of the curve portrait associated with the equivalent systems in \eqref{eq:nonind_zero} (left) and in \eqref{eq:ind_zero}~\!(right). We mark with green crosses the equilibrium points along the lines $x=0$ and $y=0$ represented in magenta. }
    \label{fig:indeterminacy}
\end{figure}

A point of indeterminacy in the above system is called a \textit{base point} if the vector field at this point is degenerate, i.e.\ if infinitely many integral curves of the field coalesce at this point. For this to be the case, at a point of indeterminacy, when re-written in any other chart of the surface in which it is contained, the vector field must also be indeterminate. The base points are identified borrowing ideas from the theory of dynamical systems. In particular, we look for equilibrium points in the ``orbitally equivalent'' system obtained by multiplying the right hand sides of the system~\!\eqref{eq:dyn_sys} by an auxiliary function\footnote{In the theory of dynamical systems, the function $g(x,y)$ should be a strictly positive smooth function. This is indeed to preserve the direction of time and not create new additional equilibrium points. Here, we relax this assumptions, since we are interested only in the geometry of the phase portrait. } $g(x,y)$ given by the least common multiple between $Q_1$ and $Q_2$ in~\!\eqref{eq:dyn_sys}. In figure~\!\ref{fig:indeterminacy} we depict the realisation of the corresponding phase portrait for the equivalent systems built for two different examples of typical systems that we encounter at this stage. We will focus on equilibrium points on the axis, as we are looking for points of indeterminacy in the original system. On the left of figure \ref{fig:indeterminacy} the vector field and its equivalent is given by the right hand sides of the following equations
\begin{equation} \label{eq:nonind_zero}
    \begin{cases}
        ~\!x' = \dfrac{x^2+y+x\,y}{x} \\[2ex]
        ~\!y'= \dfrac{1+x\,y-y^2}{x^2} 
    \end{cases} \quad \simeq \quad \begin{cases}
        ~\!x' = x^3+x\,y+x^2\,y \\[3ex]
        ~\!y'= 1+x\,y-y^2
    \end{cases}\,,
\end{equation}
while on the right of figure \ref{fig:indeterminacy} the vector field is given by the right hand sides of 
\begin{equation} \label{eq:ind_zero}
    \begin{cases}
        ~\!x' = \dfrac{x^2+y+x\,y}{x\,y^2} \\[2ex]
        ~\!y'= \dfrac{1+x\,y-y^2}{x^2} 
    \end{cases} \quad \simeq \quad \begin{cases}
        ~\!x' = x^3+x\,y+x^2\,y \\[3ex]
        ~\!y'= y^2+x\,y^3-y^4
    \end{cases}\,.
\end{equation}
In the original systems reported in \eqref{eq:nonind_zero} and \eqref{eq:ind_zero}, it seems that the points $(0,0),(0,\pm1)$ are points of indeterminacy for both systems. By looking at {the} curve portrait for the equivalent systems we can establish instead that the point $(0,0)$ is a base point for the system in \eqref{eq:ind_zero} only. { It is worth mentioning that to decide whether a point of indeterminacy is a base point for the system can also be addressed following the classical Painlev\'e test, establishing the existence of a one-parameter family of solutions locally written as a Taylor series passing through the base point. }

\begin{figure} 
  \centering
      \includegraphics[width=.65\textwidth]{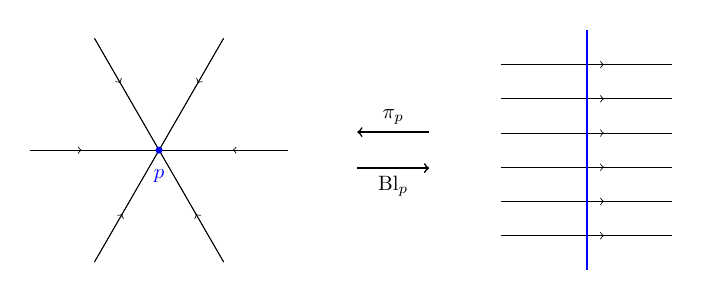}
        \caption{ With the blow-up transformation the point $p$ is replaced by the line in blue, a curve $\mathbb{CP}^1$. All the directions of the flow lines converging at the point $p$ are represented by distinct points on the blue curve after the transformation.  }
        \label{fig:directions}
    \end{figure}

The next step in the procedure of constructing the space of initial conditions is to remove all base points via the process of \textit{blowing up} the surface at such points. A blow-up is a bi-rational transformation that, in a certain sense, separates out the directions through the point in question and thus regularises the vector field at the point, introducing an \textit{exceptional line} in the process. Each point on the exceptional line (equivalent to a $\mathbb{CP}^1$) corresponds to a direction through the base point before the blow-up as schematically depicted in figure~\!\ref{fig:directions}. The definition of the space $\mathbb{C}^2$ blown up at a point $p = (a,b)$ is
\begin{equation}
  \text{Bl}_p(\mathbb{C}^2) = \{ (u_i,v_i) \times [w_0:w_1] \in \mathbb{C}^2 \times \mathbb{CP}^1: (u_i-a) w_0 = (v_i-b) w_1 \}.
\end{equation}

In coordinates, the blow-up is implemented via a bi-rational change of variables between a chart $(u_i,v_i)$ containing the base point $p=(a,b)$ and two new charts $(u_j,v_j)$, $(U_j,V_j)$, according to 
\begin{equation}
\begin{cases}
    u_i = u_j+a = U_j V_j+a  \\[1ex] 
    v_i = u_j\, v_j+b = V_j+b
\end{cases}.
\end{equation}
Here, $U_j = {w_0}/{w_1}$ covers the part of $\mathbb{CP}^1$ where $w_1 \neq 0$, while $v_j = {w_1}/{w_0}$ covers the part where $w_0 \neq 0$. We define the projection onto the first component,
\begin{equation}
    \pi_p : \text{Bl}_p(\mathbb{C}^2) \to \mathbb{C}^2, \quad (u_i,v_i) \times [w_0,w_1] \mapsto (u_i,v_i).
\end{equation}
Away from the point $p$, the blow-up is a one-to-one map,
\begin{equation}
    \text{Bl}_p(\mathbb{C}^2) \setminus \pi_p^{-1}(p) \quad \longleftrightarrow \quad \mathbb{C}^2 \setminus \{ p\}.
\end{equation}
The set $E = \pi_p^{-1}(p)$ is the exceptional curve introduced by the blow-up which in coordinates is given by {}
\begin{equation}
    E = \{ u_j = 0 \} \cup \{ V_j = 0 \}.
\end{equation}
An exceptional curve arising from the blow-up of a base point of a vector field is said to be \textit{inaccessible} if no solution of the differential system (rewritten in the according coordinates) passes through this line (i.e.\ the integral curve of the vector field do not intersect this line in a transversal way). Whether or not a line in the phase space of a Hamiltonian system is accessible can be shown by the introduction of a certain auxiliary function (approximate first integral), which, combined with a simple integral estimate, shows that this line cannot be reached by the solution, as it will be presented in section~\!\ref{sec:auxiliary_function}. 

Sometimes it will be necessary to blow down a line in the intersection diagram, in particular to compare two diagrams. This is essentially the inverse of a blow-up, but can only {be} applied to certain lines. Since the exceptional curve arising from a blow-up is a curve of self-intersection number $-1$ (see section~\!\ref{picard}), only such curves can be blown down. Furthermore, these have to be inaccessible by the vector field. Note that, in general, we do not blow down lines that have initially arisen as exceptional curves from a blow-up, as this would simply reverse this process. 

\subsection{Auxiliary function}\label{sec:auxiliary_function}
In this section we will make the observation from the previous section, that certain exceptional {curves} are inaccessible for the flow of the vector field, more rigorous. We will show that the exceptional curves introduced by each intermediate blow-up are inaccessible for the solutions of the Hamiltonian system, away from any possible base points on this curve. By intermediate blow-up we mean any blow-up that is not the final one in any cascade of blow-ups. This is done by introducing an auxiliary function $W$, that satisfies the conditions of the following lemma~\!\cite{KeckerFilipuk}.

\begin{lemma}
	\label{log_bounded}
	Suppose a function $W(z)$ is defined in the neighbourhood $U$ of a point $z_\ast$ such that the logarithmic derivative $(\log W)' = {W'}/\,{W}$ is bounded, say by $K$, on $U$. Let $\gamma \subset U$ be a finite-length curve from some point $z_0$ where $W(z_0)$ is finite and non-zero, ending in $z_\ast$. By the estimate,
	\begin{equation*}
		|\log W(z_\ast)| \leq |\log W(z_0)| + \int_{\gamma} \left| \frac{W'}{W} \right| ds \leq |\log W(z_0)| + K \,||\gamma||,
	\end{equation*} 
	$\log W(z_\ast)$, and hence $W(z_\ast)$, is bounded.
\end{lemma}

If we can find a function $W(z) = W(x(z),y(z);z)$ so that $W$ is infinite on each exceptional curve introduced by the subsequent blow-ups, while its logarithmic derivative is bounded on these lines, we can use this function in Lemma \ref{log_bounded} to show that these lines are inaccessible for the solution. Namely, for any path $\gamma \in \mathbb{C}$ with end point $z_\ast$ such that $(u_i(z_\ast),v_i(z_\ast))$ is a point on the line at infinity of $\mathbb{CP}^2$ or an exceptional curve after an intermediate blow-up, other than a base point, Lemma \ref{log_bounded} shows that $W$ must be bounded, a contradiction. This essentially shows that the exceptional curves are inaccessible for the flow of the vector field, apart from the final exceptional curve of a cascade of blow-ups, where $W$ must be finite. In particular at any movable singularity of the original system, {the solution passes transversally through a point on the exceptional curve from a final blow-up in some cascade.} In the case of the Painlev\'e equations, one finds a regular initial value problem on final exceptional curves, with analytic solution in a neighbourhood of any point on this curve, showing that the solution of the original system has a pole there. In the quasi-Painlev\'e case, the system of equations after a final blow-up becomes regular only after exchanging the role of the dependent variable, meaning we obtain a regular initial value problem for $z$ and e.g.\ {$v_{\text{f}}$ as a function of $u_{\text{f}}$ in the final chart $(u_{\text{f}},v_{\text{f}})$, yielding} a well-defined solution there. This last step is not necessary for equations with the Painlev\'e property, as in this case the system of equations is already regular and we find analytic solutions on the exceptional curve that translate into poles of the solutions, while in the quasi-Painlev\'e case we have algebraic poles.

We note that the function $W$ acts as an `approximate integral of motion', in the sense that it remains bounded at movable singularities, although it is not a constant. Such auxiliary functions were also introduced in various proofs of the (quasi-)Painlev\'e property for (quasi-)Painlev\'e equations to show that certain quantities are controlled.  
If the Hamiltonian $H$ is $z$-dependent, $H$ itself is not an integral of motion. However, for the Hamiltonian system of quasi-Painlev\'e type considered in this article, we will construct a function $W$ from the Hamiltonian by adding certain correction terms to ensure that $W$ remains bounded at any movable singularity. One then needs to check that $W$, in the transformed coordinates after each blow-up satisfies the requirements of Lemma \ref{log_bounded}. This is best done using computer algebra as the expressions for these functions can become rather large. We will demonstrate how this is done for a particular case in Appendix \ref{app:auxiliary_function}.

\subsection{Hamiltonian nature of the equations}
\label{symplectic_nature}
We start by considering a Hamiltonian system given by the function $H\big(x(z),y(z);z\big)$ and
\begin{equation*}
        y'= \dfrac{\partial H}{\partial x}  \,, \qquad 
        x'= -\dfrac{\partial H}{\partial y}\,, 
\end{equation*}
with respect to the $2$-form $\omega = dy \wedge dx$. In the subsequent changes of variables in the process of extending {the phase space by blowing up the phase space at base points, for most blow-ups introduced} we find new systems of equations in the new variables $(x,y)\mapsto (u,v)=(u(z),v(z))$ of the form 
\begin{equation}
\label{rational_equations}
f(u,v;z)\, u' = \frac{\partial H(u,v;z)}{\partial v}\,, \quad f(u,v;z)\,v' = - \frac{\partial H(u,v;z)}{\partial u}\,, 
\end{equation}
{ i.e.\ the equations remain in Hamiltonian form but with respect to a more general $2$-form $f(u,v;z) \,du \wedge dv$. 
However, the Hamiltonian form of a system of equations is not preserved under general coordinate transformations, in particular not for all $z$-dependent blow-ups, but only very special ones. }Suppose we apply a change of variables $u=u(w(z),t(z);z)$, $v=v(w(z),t(z);z)$. The system of equations becomes
\begin{equation}
\begin{aligned}
\frac{\partial u}{\partial w}\, w' + \frac{\partial u}{\partial t}\, t' + \frac{\partial u}{\partial z} &= \left. f(u,v;z)^{-1} \left( \frac{\partial H}{\partial v} \right)\right|_{u=u(w,t;z),v=v(w,t;z)} \,,  \\[2ex]
\frac{\partial v}{\partial w}\, w' + \frac{\partial v}{\partial t}\, t' + \frac{\partial v}{\partial z} &= \left. - f(u,v;z)^{-1} \left( \frac{\partial H}{\partial u} \right)\right|_{u=u(w,t;z),v=v(w,t;z)}\,. 
\end{aligned}
\end{equation}
We can solve these as equations for $(w,t)$ by
\begin{equation} \label{eq:intermediate_sympl}
\begin{aligned}
J(w,t;z)\, w' &= f(u,v;z)^{-1} \underbrace{\left( \frac{\partial H}{\partial v} - f(u,v;z) \frac{\partial u}{\partial z} \right)}_{(a)} \frac{\partial v}{\partial t} + \underbrace{\left( \frac{\partial H}{\partial u} + f(u,v;z) \frac{\partial v}{\partial z} \right)}_{(b)} \frac{\partial u}{\partial t}\,, \\[1.5ex]
-J(w,t;z)\, t' &= f(u,v;z)^{-1} \underbrace{\left( \frac{\partial H}{\partial v} - f(u,v;z) \frac{\partial u}{\partial z} \right)}_{(a)} \frac{\partial v}{\partial w} + \underbrace{\left( \frac{\partial H}{\partial u} + f(u,v;z) \frac{\partial v}{\partial z} \right)}_{(b)} \frac{\partial u}{\partial w},
\end{aligned}
\end{equation}
where 
\begin{equation}
    J(w,t;z) = \left( \frac{\partial u}{\partial w} \,\frac{\partial v}{\partial t} - \frac{\partial u}{\partial t} \,\frac{\partial v}{\partial w} \right)
\end{equation}
is the Jacobian of the transformation. { In the following we require the transformed equations to remain in Hamiltonian form, but with respect to a more general $2$-form $\omega$.
For this to be the case, the terms $(a)$ and $(b)$ in~\!\eqref{eq:intermediate_sympl} must derive from a new Hamiltonian, $K=H+h$, where $h$ satisfies the equations}
\begin{equation}
\label{h_compatible}
\frac{\partial h}{\partial v} = - f(u,v;z) \,\frac{\partial u}{\partial z}\,, \qquad \frac{\partial h}{\partial w} = f(u,v;z) \, \frac{\partial v}{\partial z}\,,
\end{equation}
that must be compatible. We can thus compute the correction $h$ to the Hamiltonian by integrating either of these equations, up to an overall function that is only dependent on $z$. Then, with the new Hamiltonian $\tilde{K}(w,t;z) = K\big(u(w,t,z),v(w,t,z);z\big)$ the equations take on the form 
\begin{equation}
F(w,t;z) \, w' = \frac{\partial \tilde{K}}{\partial t}\,, \qquad F(w,t;z) \, t' = -\frac{\partial \tilde{K}}{\partial w}\,,
\end{equation}
where $F(u,t;z) = J(w,t;z) \,\tilde{f}(w,t;z)$ and $\tilde{f}(w,t;z) = f\big(u(w,t,z),v(w,t,z);z\big)$. 
The compatibility of the equations~\!\eqref{h_compatible} severely restricts the allowed coordinate transformations for the structure of the equations { to remain in Hamiltonian form} (e.g., if the change of coordinates is not $z$-dependent). We will see below that certain blow-ups form another case where the equations remain compatible. 

One could also consider the change of variables from the point of view of the new Hamiltonian, $\tilde{K}$. The correction term to $\tilde{K}$ must satisfy
\begin{equation}
\frac{\partial \tilde{h}}{\partial w} = \tilde{f}(w,t;z) \left( \frac{\partial u}{\partial w} \frac{\partial v}{\partial z} - \frac{\partial v}{\partial w} \frac{\partial u}{ \partial z} \right)\,, \qquad \frac{\partial \tilde{h}}{\partial t} = \tilde{f}(w,t;z) \left( \frac{\partial u}{\partial t} \frac{\partial v}{\partial z} - \frac{\partial v}{\partial t} \frac{\partial u}{ \partial z} \right)\,.
\end{equation}
We note that these equations are compatible if both $f(u,v)$ and $F(w,t)$ are not explicitly $z$-dependent. 

By certain intermediate transformations one can achieve that {the $2$-form always retains the relatively simple form $\omega = u^k du \wedge dv$ for some integer power $k$. This is useful if one wants to construct a global atlas of the space of initial conditions}, with a Hamiltonian that remains polynomial in every chart, as is e.g.\ known for the Painlev\'e equations and as was also done in~\!\cite{FilipukStokes} for some equations of quasi-Painlev\'e type. However in this article we are not applying these intermediate transformations as we are just considering the blow-up process as an algorithmic way of constructing spaces of initial conditions, {and hence we do not explicitly construct a global atlas for each system.}

\subsection{Picard lattice and intersection form}
\label{picard}
For the notation introduced in this section we refer to standard introductory books on Algebraic Geometry such as~\!\cite{Shafarevich1}.
We introduce the group of divisors on a rational surface $X$. A divisor on a surface is a formal (finite) sum of { curves} (irreducible codimension $1$ subvarieties), 
\begin{equation}
    D = \sum_{i \in I} n_i \, L_i\,,
\end{equation}
where the $n_i$ are integers and $I$ a finite index set. The group of all divisors is denoted $\text{Div}(X)$. Two divisors are linearly equivalent if they differ by a principal divisor, i.e.\ a divisor given by the formal sum of the zero and pole sets (including multiplicity) of a meromorphic function $f$ on the surface, denoted
\begin{equation}
    (\,f\,) = \sum_{i} n_i \, Z_i - \sum_{j} m_j \, P_j\,,
\end{equation}
{
that is, $D_1 \sim D_2 \Leftrightarrow D_1 - D_2 = (\,f\,)$. 
Under linear equivalence, we obtain the divisor class group}
\begin{equation*}
    \text{Cl}(X) = \text{Div}(X) /\sim\,.
\end{equation*}
For example, we have 
\begin{equation}
\text{Cl}(\mathbb{P}^2) = \mathbb{Z} \,\mathcal{H}\,, \hspace{5ex} \text{Cl}(\mathbb{P}^1 \times \mathbb{P}^1) = \mathbb{Z} \,\mathcal{H}_1 \oplus \mathbb{Z} \,\mathcal{H}_2\,,
\end{equation}
where $\mathcal{H}$ is the hyperplane divisor class in $\mathbb{P}^2$ and $\mathcal{H}_1,\mathcal{H}_2$ are the divisor classes of $\{p\} \times \mathbb{P}^1$ and $\mathbb{P}^1 \times \{q\}$, respectively.
The group of classes of divisors of the surface $X$ with respect to  linear equivalence is the Picard group\footnote{The Picard group is the group of line bundles on $X$, $\text{Pic}(X)=H^1(X,\mathcal{O}^*)$. For a nonsingular rational surface, a standard fact in algebraic geometry~\cite{Shafarevich2} is $\text{Pic}(X)=H^1(X,\mathcal{O}^*)\simeq H^2(X,\mathbb{Z})$. This allows us to consider the form in~\!\eqref{eq:pic_example}. }, $\text{Pic}(X)$. { Here and in the following, we denote by $\mathcal{D}$ the divisor class, i.e.\ the element of the Picard group, for a divisor $D$.}

We introduce the intersection form $\cdot$ for divisor classes $\mathcal{D}_1, \mathcal{D}_2 \in \text{Pic}(X)$. 
For $\mathcal{D} \in \text{Pic}(X)$ we call $\mathcal{D} \cdot \mathcal{D}$ the self-intersection number of the divisor class $\mathcal{D}$. We have 
{
\begin{equation}
    \mathcal{H} \cdot \mathcal{H} = 1 \,,\quad \text{ in } \mathbb{CP}^2 \,, \hspace{8ex} \mathcal{H}_i \cdot \mathcal{H}_j = 1 - \delta_{ij} \,,\quad \text{ in } \mathbb{CP}^1 \times \mathbb{CP}^1\,. 
\end{equation} }
By the definition of a blow-up, an exceptional curve $E$ introduced in the process has self-intersection $\mathcal{E} \cdot \mathcal{E} = -1$, and 
{
\begin{equation}
    \begin{cases}
        \mathcal{E}_i \cdot \mathcal{H} = 0 \quad\forall i, \\
        \mathcal{E}_i \cdot \mathcal{E}_j = -\delta_{ij},
    \end{cases}~~ \text{ in } \mathbb{CP}^2 \,, \hspace{8ex} \begin{cases}
        \mathcal{E}_i \cdot \mathcal{H}_j = 0 \quad\forall i,j, \\
        \mathcal{E}_i \cdot \mathcal{E}_j = -\delta_{ij},
    \end{cases}~~ \text{ in } \mathbb{CP}^1 \times \mathbb{CP}^1\,. 
\end{equation} }
Furthermore, blowing up a surface $X$ at some point $P$, resulting in the blown-up surface $\tilde{X}$, we have
\begin{equation}
    \text{Pic}(\tilde{X}) = \text{Pic}(X) + \mathbb{Z} \,\mathcal{E}\,.
\end{equation}
Thus, blowing up e.g.\ $\mathbb{CP}^2$ consecutively at points $p_1,\dots,p_n$, introducing exceptional curves $E_1,\dots,E_n$ and denoting the resulting space by $X_n$, we have
\begin{equation}\label{eq:pic_example}
\text{Pic}(X_n) = \mathbb{Z} \,\mathcal{H} \oplus \mathbb{Z} \,\mathcal{E}_1 \oplus \cdots \oplus \mathbb{Z} \,\mathcal{E}_n\,. 
\end{equation}
We keep track of the blow-ups performed for a given equation by drawing the intersection diagram of the irreducible components of the inaccessible divisor arising in this process.

In this article, we start the process of computing the space of initial conditions of a system from $\mathbb{CP}^2$, which can be drawn schematically as below, the corners of the triangle being infinitely apart and the set $\{u_0 = 0\} \cup \{V_0 = 0\}$ forming the line at infinity. In terms of coordinates, we choose the following notation 
\vspace*{-4ex}
\begin{equation*}
\includegraphics[width=.29\textwidth,valign=c]{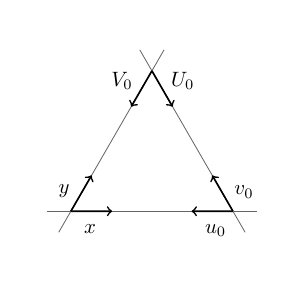} \qquad   
    \begin{aligned}
    &[\,1:x:y\,] = [\,u_0:1:v_0\,] = [\,V_0:U_0:1\,]\,, \\[.7ex]
    &u_0 = \frac{1}{x}\,, \qquad V_0 = \frac{1}{y}\,, \qquad v_0 = \frac{y}{x} = \frac{1}{U_0} \,.
\end{aligned}
\end{equation*}

\vspace*{-3ex}

A blow-up at a point $p= (u_p,v_p)$ in the chart $(u_i,v_i)$ is performed by the introduction of two new coordinate systems $(u_{i+1},v_{i+1})$, $(U_{i+1},V_{i+1})$, related to $(u_i,v_i)$ through the following transformations
\begin{equation}
    \begin{cases}
        u_{i+1} = u_i-u_p \\[1.5ex]
        v_{i+1} = \dfrac{v_i-v_p}{u_i-u_p}
    \end{cases} \,,\qquad      
    \begin{cases}
        U_{i+1} = \dfrac{u_i-u_p}{v_i-v_p} \\[1.5ex]
        V_{i+1} = v_i-v_p
    \end{cases} \,,
\end{equation}
which can easily be inverted
\begin{equation}
    \begin{cases}
    u_i = u_{i+1} + u_p = U_{i+1} V_{i+1}+u_p \\[.7ex] 
    v_i = u_{i+1} \,v_{i+1}+v_p= V_{i+1}+v_p
    \end{cases} \,.
\end{equation}
In this notation, the exceptional curve of the blow-up is identified by the set
\begin{equation}
    E = \{u_{i+1}=0\} \cup \{V_{i+1}=0\}.
\end{equation}
Usually, since we blow-up points on the exceptional curves, at least one of the coordinates of the base point is zero, i.e.\ $u_p = 0$ or $V_p = 0$. 

We note that many authors start the blow-up process from $\mathbb{CP}^1 \times \mathbb{CP}^1$, with four initial coordinate charts (e.g.\ see~\!\cite{Dzhamay2021,Dzhamay2109.06428,Takano97,Takano99}). Although the resulting diagrams are different, they contain essentially the same information. They are related by the fact that $\mathbb{CP}^1 \times \mathbb{CP}^1$ is obtained from $\mathbb{CP}^2$ by two blow-ups (at the points $(u_0,v_0)=(0\,,0)$ and $(U_0,V_0)=(0\,,0)$), giving rise to two exceptional curves $E_1$ and $E_2$, followed by a blow-down of the emerging $-1$ curve $H-E_1-E_2$.

We now outline the blow-up process in terms of divisors. For each blow-up, two situations can occur. 
\begin{enumerate}
    \item The new base point lies on the intersection of two irreducible components of the inaccessible divisor $L_1$ and $L_2$
     \vspace*{-4ex}
    \begin{equation*}
        \includegraphics[width=.7\textwidth]{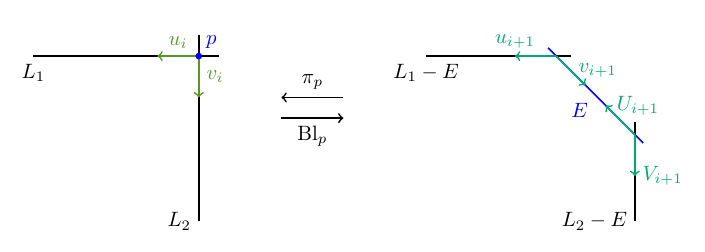}
    \end{equation*}
    
    \vspace*{-3ex}
    
    In this case, the base point is at the origin $(u_i,v_i)=(0\,,0)$ and not \textit{visible} in other charts. After the blow-up, the irreducible divisor components of the total transform $\hat{L}_1 \cup \hat{L}_2$ are $\hat{L}_1 - E$, $E$ and $\hat{L}_2 - E$, where $E$ is the exceptional curve of the blow-up and the notation $L_i-E$ stands for the strict transform of $L_i$ under the blow-up, i.e.\ 
    \begin{equation*}
        L_i-E = \overline{\pi_p^{-1}\{L_i \setminus\{p\}\}} \,.
    \end{equation*}
    In the following we will omit the hats and simply write $L_i-E$ by abuse of notation.

    \item The new base point is not an intersection point of inaccessible divisor components and therefore away from the origin of any coordinate chart 
    \vspace*{-4ex}
     \begin{equation*} 
        \includegraphics[width=.7\textwidth]{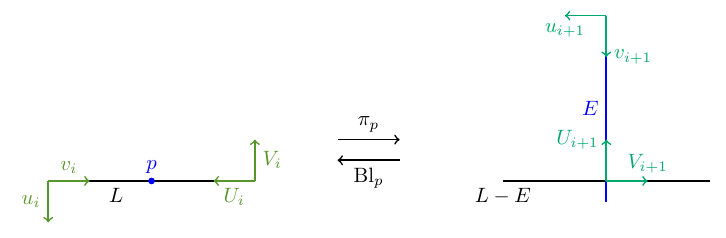}
    \end{equation*}
     
    \vspace*{-3ex}
    
    Therefore, it will be visible in both charts 
    $$(u_i,v_i) = (0\,,a)\,, \qquad (U_i,V_i) = (1/a,0)\,, \quad a \in \mathbb{C}\,.$$
\end{enumerate}
Each base point of the extended Hamiltonian system in $\mathbb{CP}^2$ gives rise to a \textit{cascade} of blow-ups. Sometimes a cascade can split into several branches, namely if several new base points appear on the exceptional curve after a blow-up. A cascade of blow-ups finishes if on the exceptional curve after a blow-up both vector fields $(u_\text{f}',v_\text{f}')$, $(U_\text{f}',V_\text{f}')$ form a regular initial value problem after a change of independent variables, i.e.\ no new base points appear. We call the corresponding charts $(u_\text{f},v_\text{f})$, $(U_\text{f},V_\text{f})$ the \textit{final charts} of the cascade. { In all cases that we consider in this paper, each cascade of blow-ups finishes after a finite number of blow-ups. Although a resolution of base points is possible in a finite number of steps, we do not know of an a priori criterion for an equation to decide how many blow-ups are necessary.}

For a given equation we obtain the blow-up diagram which schematically depicts the irreducible components of the exceptional divisor. Each line in this diagram is assigned its self-intersection number. In a further abstraction we represent each line by a dot, where we connect dots whose corresponding lines have a non-zero intersection product. In the case of Painlev\'e equations, these are (extended) Dynkin diagrams (all irreducible components have self-intersection $-2$), whereas for the quasi-Painlev\'e equations studied in section~\ref{sec:quasi_Painleve}, we have at least one curve of self-intersection $-3$.

From the form of the system of equations in the final charts of a cascade of blow-ups we can deduce the possible types of movable singularity that the system can develop when passing through the final exceptional curve.

\section{Hamiltonian systems with the Painlev\'e property}

Hamiltonian forms for all six Painlev\'e equations appear already in an article as far back as 1923 by Malmquist \cite{Malmquist}. 
Much later, Okamoto \cite{okamoto1,okamoto2,okamoto3,okamoto4} studied all Painlev\'e equations in their Hamiltonian form and their (auto-)B\"acklund symmetries as affine Weyl groups. These are discrete symmetries of the equations that map solutions of one system to the same equations but with different parameters. That is, they act on the space of parameters of the equation. 
In a subsequent paper~\cite{Okamoto1979}, Okamoto introduced the space of initial conditions for each of the Painlev\'e Hamiltonian systems, compactifying the phase space to some Hirzebruch surface, applying $8$ blow-ups and removing the inaccessible divisors from the space. He showed that the solution curves of the Painlev\'e flows uniformly foliate these spaces, with leaves {transversal} to each fibre. 
Yet another couple of decades later, H. Sakai published the seminal paper~\cite{Sakai2001}, classifying discrete and differential Painlev\'e systems in terms of affine Weyl groups belonging to their surface type and symmetry type, thus establishing a kind of dual relationship between these two. 
The Painlev\'e identification problem is to establish, given an equation or system with the Painlev\'e property, to which of the equations in the list of Painlev\'e equations it is related by bi-rational transformations. 
In order to explain the geometric approach and introduce further notation, we start by reviewing the case of Hamiltonian systems related to the Painlev\'e-II equation. There are two Hamiltonian forms for this equation usually considered in the literature, but below we have indeed listed four Hamiltonians giving rise to this equation, establishing pairwise bi-rational transformations between any of them through the geometric approach.
Rather than starting from Hirzebruch surfaces for compact phase spaces as Okamoto did, here and throughout the article we start the process of removing points of indeterminacy in the system from $\mathbb{CP}^2$. Although this leads to the transformations between the resulting systems after final blow-ups to be non-symplectic, it has the advantage for us that the process of obtaining the space of initial conditions is completely algorithmic. The process finishes once one has obtained a regular initial value problem in each chart where the equations are defined, possibly after a change of dependent variables.

\subsection{Hamiltonian systems for Painlev\'e-II}
\label{P2review}
In this section we review the construction of Okamoto's space of initial conditions for several systems relating to the Painlev\'e-II equation, starting from the compact phase space $\mathbb{CP}^2$. For each system, we obtain a blow-up diagram, showing the configuration of exceptional curves arising from the blow-up process to regularise the system. By comparing diagrams for two equivalent systems, we demonstrate how to obtain a bi-rational change of variables between these. This is the geometric approach to identifying Hamiltonian systems. The aim is to then apply this approach to a number of systems of quasi-Painlev\'e type in section \ref{sec:quasi_Painleve}.

We consider the Painlev\'e-II equation in Hamiltonian form. Namely, we have the four Hamiltonians $H^{\text{P}_{\text{I\!I}}}_k(x_k(z),y_k(z);z)$, $k \in \{1,2,3,4\}$, each of which gives rise to a different system of equations, all of which however reduce to the Painlev\'e-II equation in $y_k$, under elimination of $x_k$. These are
\begin{align} 
    H_1^{\text{P}_{\text{I\!I}}}(x_1,y_1;z) & = \frac{1}{2} \,x_1^2 - \left(y_1^2 + \frac{z}{2} \right) x_1 - \left( \alpha + \frac{1}{2} \right) y_1 \,,  \label{eq:Ham_P2_1} \\[.5ex]
    H_2^{\text{P}_{\text{I\!I}}}(x_2,y_2;z) & = \frac{1}{2} \,x_2^2 - y_2^2 \, x_2 - \frac{z}{2}\, y_2^2 - \alpha \,y_2 \,, \label{eq:Ham_P2_2} \\[.7ex]
    H_3^{\text{P}_{\text{I\!I}}}(x_3,y_3;z) & = \frac{1}{2}\,x_3^2 - \frac{1}{2} \,y_3^4 - \frac{z}{2} \, y_3^2 - \alpha \,y_3 \label{eq:Ham_P2_3} \,, \\[.7ex]
    H_4^{\text{P}_{\text{I\!I}}}(x_4,y_4;z) & = \frac{3}{4}\, x_4^2 - 2 \,x_4\, y_4^2 + y_4^4   - \frac{z}{4}\, x_4 - \frac{1+4\,\alpha}{6}\, y_4 \, \label{eq:Ham_P2_4} .
\end{align}
\begin{figure}[t]
\centering
\includegraphics[width=.9\textwidth]{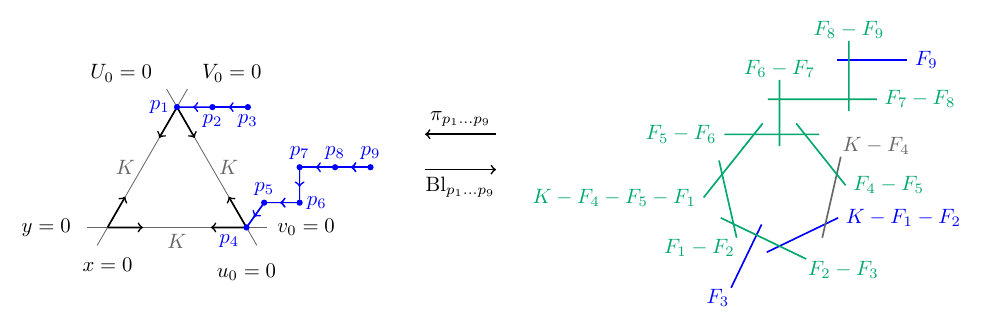} 
\caption{Blow-up cascades for $H_1^{\text{P}_{\text{I\!I}}}$ and space of initial conditions. On the right, in green the curves with self-intersection $-2$, in blue with self-intersection $-1$, in gray with self-intersection $\ge 0$.}
    \label{fig:Ham_P2_1_1st_blowup}
\end{figure}
\begin{figure}[t]
    \centering
   \includegraphics[width=.9\textwidth]{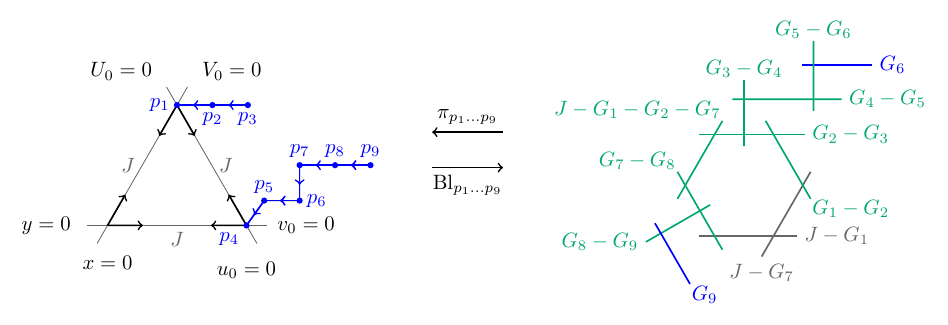}
    \caption{Blow-up cascades for $H_2^{\text{P}_{\text{I\!I}}}$ and space of initial conditions. On the right, in green the curves with self-intersection $-2$, in blue with self-intersection $-1$, in gray with self-intersection $\ge 0$.}
    \label{fig:Ham_P2_2}
\end{figure}
The coordinate transformations between the Hamiltonian systems obtained from $H_k^{\text{P}_{\text{I\!I}}}$, $k = 1,2,3$ are easy enough to be guessed directly. Namely,
\vspace*{-1ex}

\begin{equation}
\begin{cases}   x_2  = x_1 - \dfrac{z}{2}\\[.7ex]  
     y_2 = y_1 
     \end{cases}\,, \qquad \qquad 
     \begin{cases}
        x_3 = x_2 - y_2^2\\[.7ex]   
        y_3 = y_2 
        \end{cases} \,, \qquad \qquad 
     \begin{cases}
        x_4 = \dfrac{2}{3}\,x_3 + \dfrac{4}{3}\,y_3^2+\dfrac{z}{6}   \\[.7ex]
        y_4 = y_3
        \end{cases} \,. 
\end{equation}

However, we demonstrate here how these transformations can also be found using the geometric approach to the Painlev\'e equations, by compactifying each of these systems on $\mathbb{CP}^2$, applying the blow-up procedure and identifying irreducible components of inaccessible divisors. For Painlev\'e-II this is already done in~\!\cite{Dzhamay2021}, however starting from $\mathbb{CP}^1 \times \mathbb{CP}^1$. This identification process is also exhibited excellently in~\!\cite{Dzhamay2109.06428} for the case of the Painlev\'e-IV equation. Our aim is to apply this identification method in the next section to a number of Hamiltonian systems of quasi-Painlev\'e type, essentially to obtain a classification of systems with a certain singularity structure.

We construct the surface diagram for $H_1^{\text{P}_{\text{I\!I}}}$ in~\eqref{eq:Ham_P2_1} starting from the identification of the base points to blow-up, represented in figure~\ref{fig:Ham_P2_1_1st_blowup}.
Initially, there are two base points in the system of equations compactified on $\mathbb{CP}^2$ located at 
\begin{equation}
p_1\colon (U_0,V_0)=(0\,,0) \,, \qquad p_4\colon (u_0,v_0) =(0\,,0) \,.
\end{equation}
The two cascades of blow-ups obtained from these are
\begin{align}
    p_1 \,\, &\longleftarrow\,\, p_2\colon (U_1,V_1) = (0\,,0)\,\, \longleftarrow\,\, p_3\colon (U_2,V_2)=\!\Big(-\alpha-\frac{1}{2}\,,0\Big)\,, \\[1ex]
\begin{split}
    p_4 \,\, &\longleftarrow\,\, p_5\colon(U_4,V_4) = (0\,,0)  \,\, \longleftarrow\,\, p_6\colon \!\Big(u_5,v_5\Big) = \!\Big(0\,,\frac{1}{2}\Big) \,\, \longleftarrow\,\, p_7\colon(u_6,v_6)= (0\,,0) \\[.7ex] 
    &\longleftarrow\,\, p_8\colon \!\Big( u_7,v_7\Big)= \!\Big( 0\,, -\frac{z}{4}\Big)  \,\, \longleftarrow\,\, p_9\colon\!\Big( u_8, v_8\Big) = \!\Big(0\,,\frac{1}{8}\,(1-2\alpha) \Big), 
\end{split}
\end{align}
after performing all the blow-ups the space of initial conditions is schematised in figure~\!\ref{fig:Ham_P2_1_1st_blowup}. 
    Since the case for $H_2^{\text{P}_{\text{I\!I}}}$ is very similar to the one for $H_1^{\text{P}_{\text{I\!I}}}$ we do not explicitly write out the blow-up cascades here, but only draw its surface diagram (figure~\!\ref{fig:Ham_P2_2}). 

\begin{figure}[t]
    \centering
 \includegraphics[width=.95\textwidth]{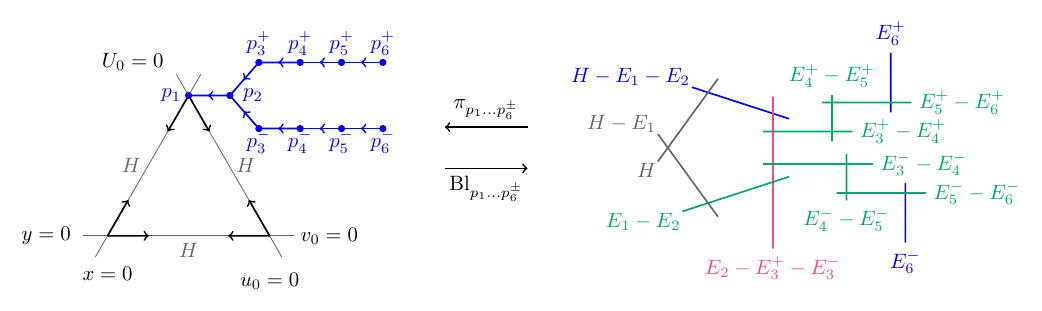}
    \caption{Blow-up cascades for $H_3^{\text{P}_{\text{I\!I}}}$ and $H_4^{\text{P}_{\text{I\!I}}}$ and relative space of initial conditions. On the right, in magenta the curves with self-intersection $-3$, in green with self-intersection $-2$, in blue with self-intersection $-1$, in gray with self-intersection $\ge 0$.}
    \label{fig:Ham_P2_3}
\end{figure}

    The case of $H_3^{\text{P}_{\text{I\!I}}}$ is somewhat different and we give the cascades of blow-ups in the following. After the second blow-up there are two new base points,  cascade of blow-ups, $q_3^{+}\colon (u_2^{+},v_2^{+}) = (0\,,1)$ and $q_3^{-}\colon (u_2^{-},v_2^{-}) = (0\,,-1)$, that is the cascade of blow-ups splits into two branches:
\begin{align} 
\begin{split} \label{eq:PII_H3_casc1}
    q_1&\colon  (u_0,v_0) = (0\,,0)\,\, \longleftarrow\,\, q_2\colon (U_1,V_1) = (0\,,0)\,\, \longleftarrow\,\, q_3^{+}\colon (u_2^{+},v_2^{+}) = (0\,,1) \\[.7ex] & \longleftarrow\,\, q_4^{+}\colon (u_3^{+},v_3^{+}) = (0\,,0)\,\, \longleftarrow\,\, q_5^{+}\colon (u_4^{+},v_4^{+}) = \!\Big( 0\,, -\frac{z}{2} \Big)\,\, \longleftarrow\,\, q_6^{+}\colon (u_5^{+},v_5^{+}) = \!\Big( 0\,, \frac{1}{2} - \alpha \Big),
    \end{split} \\[2ex]
 \begin{split} \label{eq:PII_H3_casc2}
    q_2 &\,\, \longleftarrow\,\, q_3^{-}\colon (u_2^{-},v_2^{-}) = (0\,,-1)\,\, \longleftarrow\,\, q_4^{-}\colon (u_3^{-},v_3^{-}) = (0\,,0)\,\, \longleftarrow\,\, q_5^{-}\colon (u_4^{-},v_4^{-}) = \!\Big( 0\,, \frac{z}{2} \Big) \\ &\,\, \longleftarrow\,\, q_6^{-}\colon (u_5^{-},v_5^{-}) = \!\Big( 0\,, \frac{1}{2} + \alpha \Big).
    \end{split} 
\end{align} 
The corresponding surface diagram is depicted in figure~\!\ref{fig:Ham_P2_3} and it is the same diagram obtained for $H_4^{\text{P}_{\text{I\!I}}}$ in~\!\eqref{eq:Ham_P2_4}. In the latter case the cascade is slightly different from the first one, for the two branches starting from two different points, as
\begin{align}
\begin{split}
    q_1&\colon  (u_0,v_0) = (0\,,0) \,\, \longleftarrow \,\, q_2\colon (u_1,v_1) = (0\,,0) \,\, \longleftarrow \,\, q_3^{+}\colon (u_2^{+},v_2^{+}) = \!\Big(0\,,\frac{1}{2}\Big) \\[.7ex] 
    & \hspace*{-2.5ex} \longleftarrow \,\, q_4^{+}\colon (u_3^{+},v_3^{+}) = (0\,,0) \,\, \longleftarrow \,\, q_5^{+}\colon (u_4^{+},v_4^{+}) = \!\Big( 0\,,-\frac{z}{8} \Big) \,\, \longleftarrow \,\, q_6^{+}\colon (u_5^{+},v_5^{+}) = \!\Big( 0\,,\frac{1}{12}(1-2\,\alpha)\Big)
\end{split} \\[1ex]
\begin{split}
    q_2&\colon (u_1,v_1) = (0\,,0) \,\, \longleftarrow \,\, q_3^{-}\colon (u_2^{-},v_2^{-}) = \!\Big(0\,,\frac{3}{2}\Big) \,\, \longleftarrow \,\, q_4^{-}\colon (u_3^{-},v_3^{-}) = (0\,,0) \\[.7ex] 
    &  \hspace*{-2.5ex} \longleftarrow \,\, q_5^{-}\colon (u_4^{-},v_4^{-}) = \!\Big( 0\,, \frac{3}{8}\,z \Big) \,\, \longleftarrow \,\, q_6^{-}\colon (u_5^{-},v_5^{-}) = \!\Big( 0\,, \frac{3}{4}(1+2\,\alpha)\Big). 
\end{split}
\end{align} 

The surface diagram for $H_3^{\text{P}_{\text{I\!I}}}$ needs to be transformed into its minimal form. To achieve this, we proceed by blowing down the $-1$ curve $H-E_1-E_2$, so that the $-3$ curve $E_2 - E_3^+ -E_3^-$ is transformed into a $-2$ curve. 

\subsection{Comparison of two Hamiltonian systems for \texorpdfstring{$\text{P}_{\text{II}}$}{pII}}
In order to compare the different systems, we refer to the Dynkin diagram ${E}_7^{(1)}$, which emerges by considering the Picard lattice associated with the surface in its minimal form. For $H_1^{\text{P}_{\text{I\!I}}}$ we have 
\begin{equation*}
    \text{Pic}(X) = \text{Span}_{\mathbb{Z}}
    \,\{ \mathcal{K} ,\mathcal{F}_1 \,,\dots, \mathcal{F}_9 \}
\end{equation*}
\begin{equation*}
    \begin{tikzpicture}[scale=.75,baseline={([yshift=-.5ex]current bounding box.center)}]
        \coordinate (o1) at (0,-1);
        \draw (o1) -- coordinate (M) ++(6,0);
        \draw (M) -- ++(0,1) coordinate (N); 
        \foreach \i in {1,...,7} {
        \path[draw=black,fill=white] (\i-1,-1) circle (2pt) node[below=.1] {$\delta_{\i}$};
        } 
        \path[draw=black,fill=white] (N) circle (2pt) node[above=.1] {$\delta_8$};
    \end{tikzpicture} \hspace{15ex}  
    \begin{aligned}
    &\delta_1 \colon \mathcal{F}_8 - \mathcal{F}_9 \qquad  &\delta_5 &\colon \mathcal{K} - \mathcal{F}_4 - \mathcal{F}_5 - \mathcal{F}_1 \\
    &\delta_2 \colon \mathcal{F}_7 - \mathcal{F}_8  \qquad &\delta_6 &\colon \mathcal{F}_1 - \mathcal{F}_2 \\
    &\delta_3 \colon \mathcal{F}_6 - \mathcal{F}_7 \qquad &\delta_7 &\colon \mathcal{F}_2 - \mathcal{F}_3 \\
    &\delta_4 \colon \mathcal{F}_5 - \mathcal{F}_6   \qquad &\delta_8 &\colon \mathcal{F}_4 - \mathcal{F}_5  \\
    \end{aligned}
\end{equation*}
We now compare the irreducible components of the inaccessible divisor of any two of these systems. We exemplarily choose $H_1^{\text{P}_{\text{I\!I}}}$ and $H_3^{\text{P}_{\text{I\!I}}}$ to demonstrate the process. Note that the diagram for $H_3^{\text{P}_{\text{I\!I}}}$ is not minimal in the sense that it contains an inaccessible $-1$ curve, namely the line $H-E_1-E_2$ in figure~\!\ref{fig:Ham_P2_3}, that can be removed by a blow-down. To compare the diagrams for $H_1^{\text{P}_{\text{I\!I}}}$ and $H_3^{\text{P}_{\text{I\!I}}}$, we can alternatively blow up the diagram for $H_1^{\text{P}_{\text{I\!I}}}$ at an additional point which has to be on the line $F_5-F_6$. We denote the additional exceptional curve thus obtained by $F_{10}$.

We can now compare irreducible divisor components as follows. There is some freedom in the choice made below, where we identify $\mathcal{F}_9 = \mathcal{E}_6^{-}$ and $\mathcal{F}_3 = \mathcal{E}_6^{+}$. The choice the other way round would result in a slightly different coordinate transformation related to the one we obtain below by interchanging some negative signs, as
\begin{equation}
\begin{aligned}
    \mathcal{K} - \mathcal{F}_4 - \mathcal{F}_5 - \mathcal{F}_1 & = \mathcal{E}_3^{-} - \mathcal{E}_4^{-}, & \quad \mathcal{F}_6 - \mathcal{F}_7 & = \mathcal{E}_3^{+} - \mathcal{E}_4^{+}, \\
    \mathcal{F}_1 - \mathcal{F}_2 & = \mathcal{E}_4^{-} - \mathcal{E}_5^{-}, & \quad  \mathcal{F}_7 - \mathcal{F}_8 & = \mathcal{E}_4^{+} - \mathcal{E}_5^{+}, \\
    \mathcal{F}_2 - \mathcal{F}_3 & = \mathcal{E}_5^{-} - \mathcal{E}_6^{-}, & \quad \mathcal{F}_8 - \mathcal{F}_9 & = \mathcal{E}_5^{+} - \mathcal{E}_6^{+}, \\
    \mathcal{F}_4 - \mathcal{F}_5 & = \mathcal{E}_1 - \mathcal{E}_2, & \quad \mathcal{F}_{10} & = \mathcal{H} - \mathcal{E}_1 - \mathcal{E}_2, \\ 
    & \hspace{8ex} &\hspace{8ex} \mathcal{F}_5 - \mathcal{F}_6 - \mathcal{F}_{10}, & = \mathcal{E}_2 - \mathcal{E}_3^{+} - \mathcal{E}_3^{-}. \\
\end{aligned}
\end{equation}
This can be achieved by the following identifications
\begin{equation}
\begin{aligned}
    \mathcal{F}_1&=\mathcal{E}_4^{-}, ~~~ \mathcal{F}_2 = \mathcal{E}_5^{-},~~~  \mathcal{F}_3 = \mathcal{E}_6^{-},    &\hspace{12ex} \mathcal{E}_1 &= \mathcal{K} - \mathcal{F}_5 - \mathcal{F}_{10}, \\
    \mathcal{F}_4 &= \mathcal{H} - \mathcal{E}_2 - \mathcal{E}_3^{-} , &\hspace{12ex} \mathcal{E}_2 &= \mathcal{K} - \mathcal{F}_4 - \mathcal{F}_{10}, \\
    \mathcal{F}_5 & = \mathcal{H} - \mathcal{E}_1 - \mathcal{E}_3^{-}, &\hspace{12ex} \mathcal{E}_3^{-} & = \mathcal{K} - \mathcal{F}_4 - \mathcal{F}_5, \\
    \mathcal{F}_6 & = \mathcal{E}_3^{+} , ~~~ \mathcal{F}_7  = \mathcal{E}_4^{+} &\hspace{12ex} \mathcal{E}_4^{-} & = \mathcal{F}_1 , ~~~ \mathcal{E}_5^{-} = \mathcal{F}_2, ~~~ \mathcal{E}_6^{-}  = \mathcal{F}_3, \\
    \mathcal{F}_8 & = \mathcal{E}_5^{+} , ~~~ \mathcal{F}_9  = \mathcal{E}_6^{+} &\hspace{12ex}  \mathcal{E}_3^{+} & = \mathcal{F}_6, ~~~ \mathcal{E}_4^{+}  = \mathcal{F}_7 , \\
    \mathcal{F}_{10} & = \mathcal{H} - \mathcal{E}_1 - \mathcal{E}_2, &\hspace{12ex}  \mathcal{E}_5^{+} & = \mathcal{F}_8, ~~~ \mathcal{E}_6^{+}  = \mathcal{F}_9, \\ 
    \mathcal{K} & = 2\mathcal{H} - \mathcal{E}_1 - \mathcal{E}_2 - \mathcal{E}_3^{-}. &\hspace{12ex} \mathcal{H} & = 2\mathcal{K} - \mathcal{F}_1 - \mathcal{F}_7 - \mathcal{F}_{10}. 
\end{aligned}
\end{equation}
Now, the line $\{y_1=0\}$ ($x_1$-axis) in the diagram for $H_1^{\text{P}_{\text{I\!I}}}$ is given by $\mathcal{K}-\mathcal{F}_4$, whereas the line $\{x_1=0\}$ ($y_1$-axis) is $\mathcal{K}-\mathcal{F}_1-\mathcal{F}_2$. On the other hand, the $x_3$-axis in the diagram for $H_3^{\text{P}_{\text{I\!I}}}$ is $\mathcal{H}-\mathcal{E}_1$, the $y$-axis being just $\mathcal{H}$. We have 
\begin{equation}
    \mathcal{K}-\mathcal{F}_4 = 2\mathcal{H} - \mathcal{E}_1 - \mathcal{E}_2 - \mathcal{E}_3^{-} - (\mathcal{H} - \mathcal{E}_2 - \mathcal{E}_3^{-}) = \mathcal{H}-\mathcal{E}_1,
\end{equation}
meaning that $y_1$ is related to $x_3$, $y_3$ by a linear transformation (note that the divisor class $\mathcal{H}$ contains both the lines $\{x=0\}$ and $\{y=0\}$)
\begin{equation}
    y_1 = A x_3 + B y_3 + C,
\end{equation}
the line thus defined passing through the point $q_1\colon (u_0,v_0)=(0\,,0)$. This condition gives $A=0$, so that $y_1 = B\,y_3 +C$.

Furthermore, for the line $\{x_1=0\}$ we have
\begin{equation}
    \mathcal{K}-\mathcal{F}_1-\mathcal{F}_2 = 2\,\mathcal{H} - \mathcal{E}_1 - \mathcal{E}_2 - \mathcal{E}_3^{-} - \mathcal{E}_4^{-} - \mathcal{E}_5^{-},
\end{equation}
where the divisor on the right-hand side is equivalent to a curve quadratic in $x_3$ and $y_3$, 
\begin{equation}
    a \,x_3^2 + b \,x_3\,y_3 + c \,y_3^2 + d\,x_3 + e\,y_3 + f = 0.
\end{equation}
The coefficients are fixed by the fact that the curve is passing through the points $q_1$, $q_2$, $q_3^{-}$, $q_4^{-}$ and $q_5^{-}$ in~\!\eqref{eq:PII_H3_casc1} and~\!\eqref{eq:PII_H3_casc2}, these conditions resulting in $a=b=e=0$, $c=d$ and $f = c\,z/2$, so that
\begin{equation}
    x_1 = c \left(y_3^2 + x_3 + \frac{z}{2} \right)
\end{equation}
The remaining open constants $B$, $C$ and $c$ are determined by requiring that $x_1$, $y_1$ actually satisfy the system defined by $H_1^{\text{P}_{\text{I\!I}}}$, giving $B=c=1$ and $C=0$.

This means that $x_1,y_1$ are related to $x_3,y_3$ by the transformation
\begin{equation}
    x_1 = y_3^2 + x_3 + \frac{z}{2}, \qquad y_1 = y_3,
\end{equation}
the inverse of this coordinate change being
\begin{equation}
    x_3 = x_1 - y_1^2 - \frac{z}{2}, \qquad y_3 = y_1. 
\end{equation}
One can easily check that this bi-rational change of variables is symplectic, in the sense that the form
\begin{equation}
\Omega_1 = dx_1 \wedge dy_1 - d H_1^{{\text{qsi-P}}_{\text{I\!I}}} \wedge dz
\end{equation}
is preserved, i.e.\ if we denote by $\phi\colon (x_1,y_1) \mapsto (x_3,y_3)$ the change of coordinates, we have
\begin{equation}
\Omega_1 = \phi^{\ast} (\Omega_3)\,, \qquad \Omega_3 = dx_3 \wedge dy_3 - d H_3^{{\text{qsi-P}}_{\text{I\!I}}} \wedge dz.
\end{equation}

\paragraph{Remark} In addition to the change of dependent variables, it may be necessary to make a change in the independent variable $z$ to compare two symplectic structures. E.g., had we made the choice $\mathcal{F}_9 = \mathcal{E}_6^{+}$, $\mathcal{F}_6 = \mathcal{E}_6^{-1}$ above, we would obtain a slightly different change of variables, which however is only realised as a symplectic transformation if additionally we let $z=-t$, where $z$ and $t$ are the independent variables of the systems $H_1^{\text{P}_{\text{I\!I}}}$ and $H_3^{\text{P}_{\text{I\!I}}}$, respectively
\begin{comment}
\begin{equation}
\begin{aligned}
    x_1 & = y_3^2 - x_3 -\frac{t}{2}, \quad & y_1 & = y_3, \\
    x_3 &= y_1^2 - x_1 - \frac{z}{2}, \quad & y_3 &= y_1,
\end{aligned}
\end{equation}
\end{comment}
\begin{equation}
\begin{cases}
    x_1  = y_3^2 - x_3 -\dfrac{t}{2} \\[.7ex]
    y_1  = y_3
\end{cases} \,, \qquad 
\begin{cases}
    x_3 = y_1^2 - x_1 - \dfrac{z}{2} \\[.7ex]
    y_3 = y_1
\end{cases} \,.
\end{equation}

\section{Hamiltonian systems with the quasi-Painlev\'e property}\label{sec:quasi_Painleve}
Several generalisations of the Painlev\'e property have been proposed in the literature, such as the weak Painlev\'e property \cite{GrammaticosRamani}, poly-Painlev\'e property \cite{KruskalClarkson} and algebro-Painlev\'e property \cite{KeckerHalburd2014}. In work by Shimomura~\!\cite{Shimomura2006,Shimomura2008}, classes of differential equations are introduced for which it is shown that all movable singularities that can be reached by continuation of an analytic solution along a finite-length path in the plane, are algebraic poles of a certain type. This property is now generally referred to as the quasi-Painlev\'e property (from the title of the papers~\!\cite{Shimomura2006,Shimomura2008}). 
Filipuk and Halburd~\!\cite{halburd1} consider an even more general class of second-order ordinary differential equations,
\begin{equation}
\label{quasi_2ndorder}
y''(z) = y^N + a_{N-2}(z) y^{N-2} + \cdots + a_1(z) y + a_0(z),
\end{equation}
which, under certain conditions on the coefficient functions $a_0,\dots,a_{N-2}$, possess the quasi-Painlev\'e property. In the even $N$ case, this condition is $a_{N-2}''(z) \equiv 0$ (i.e.\ $a_{N-2}$ is either a linear or a constant function). In the odd $N$ case, there are two conditions, $a_{N-2}''(z) \equiv 0$ and another, slightly more involved one that we do not write down in general. E.g., for $N=5$ the other condition is $\left( a_1(z) - a_3^2(z) \right)' \equiv 0$. The movable singularities of the equations~\!\eqref{quasi_2ndorder} are given by Puiseux expansions,
\begin{equation}
\begin{aligned}
y(z) & = \sum_{j=-2}^\infty c_j (z-z_\ast)^{j/(N-1)}, \quad N \text{ even}, \\
y(z) & = \sum_{j=-1}^\infty c_j (z-z_\ast)^{2j/(N-1)}, \quad N \text{ odd}.
\end{aligned}
\end{equation}
The equations~\!\eqref{quasi_2ndorder} can be written trivially in Hamiltonian form,
\begin{equation}
	\label{2nd-order-Ham}
H(x,y;z) = \frac{1}{2}\, x^2 - \frac{1}{N+1}\, y^{N+1} - \frac{a_{N-2}(z)}{N-1}\, y^{N+1} + \cdots + \frac{a_1(z)}{2}\, y^2 + a_0(z)\, y,
\end{equation}
where one can in principle add any function in $z$, which doesn't impact the Hamilton's equations,
\begin{equation*}
y'(z) = \frac{\partial H}{\partial x}\,, \qquad x'(z) = - \frac{\partial H}{\partial y}\,.
\end{equation*}
For $N=2,3$ these Hamiltonians reduce to the ones for the Painlev\'e equations I and II, for which the space of initial conditions was given by Okamoto. In the quasi-Painlev\'e case, for $N=4$ and $N=5$, the space of initial conditions was constructed in ~\!\cite{KeckerFilipuk}. The difference to the Painlev\'e case is that, after the final blow-up in any cascade, one needs to make a hodographic transformation, exchanging the role of the independent variable $z$ with one of the coordinate functions, in order to obtain a regular initial value problem. To obtain the spaces of initial conditions, a number of $14$ blow-ups was necessary in the case $N=4$ and $15$ blow-ups in the case $N=5$ (although the surface obtained here is not minimal, and in fact a number of $13$ blow-ups suffices when the Hamiltonian form is written instead in one of the forms below).
It was shown in \cite{FilipukStokes} that, by removing all inaccessible divisors, one can obtain a global atlas for the space of initial conditions, where the coordinate transformations between different charts are (almost) symplectic, in a way that the $2$-form can have zeros in the final blow-up charts,
\begin{equation}
\label{2form-simple}
\omega = dx \wedge dy = u\, du \wedge dv.
\end{equation}
Furthermore, this can be done in such a way that the Hamiltonian functions in all charts are polynomial in their dependent variables.

In the next section, we present a number of Hamiltonian systems for which the $2$-form has the form~\!\eqref{2form-simple}. Similar as for the case of the Painlev\'e equations, we can identify some of these systems by the analogue of their Okamoto's spaces of initial conditions, namely by finding correspondences between the irreducible components of the inaccessible divisor. While in the Painlev\'e case, the inaccessible divisor corresponds to the pole divisor of the $2$-form $\omega$, this is no longer the case for quasi-Painlev\'e equations, as $\omega$ has zeros on the final exceptional curve. E.g., in the case of square-root type singularities, $\omega$ is canonical on the second to final exceptional curve of a blow-up cascade and has a simple zero on the final exceptional curve as in~\eqref{2form-simple}. In the general case of a $k$th-root type singularity, $\omega$ has a zero of order $k-1$ on the final exceptional curve as in~\eqref{quasi_symplectic_form}, zeros of lower order on the next to final exceptional curves and is canonical on the $k$th to final exceptional curve. Furthermore we note that, while in the Painlev\'e case, by the uniqueness results in~\cite{Takano97,Takano99} one is guaranteed to find a bi-rational change of variables to a known Painlev\'e equation if the configuration of the inaccessible divisor components matches the one of that equation, such a result has not been established for the quasi-Painlev\'e case in general. 

\subsection{Hamiltonian systems of quasi-Painlev\'e-II type}
\label{quasi-PainleveII}

Apart from the (trivial) Hamiltonian~\!\eqref{2nd-order-Ham}, we find several other representations for the equation
\begin{equation}
\label{quasiPII}
y'' = 3y^5 + a_3(z)\, y^3 + a_2(z)\, y^2 + a_1(z)\, y + a_0(z),
\end{equation}
the factor of $3$ multiplying the leading-order term here chosen for convenience. This equation was shown to have the quasi-Painlev\'e property under the resonance conditions $a_3''=0$ and $(a_3^2-16a_1)'=0$. Therefore $a_3$ is at most linear, i.e.\ $a_3(z) = az+b$ and in the case $a\neq0$ we can shift and re-scale $z$ such that
\begin{equation}\label{eq:contraints_quasi_PII}
    a_3(z)=z\,, \qquad a_1(z) = \frac{z^2}{16} + c \,, 
\end{equation}
with $c \in \mathbb{C}$ a constant. We remark that in the case $a=0$,~\!\eqref{quasiPII} would be a hyper-elliptic equation which can be integrated by quadrature.
Similar to the case of the Painlev\'e-II equation, we can give the following three representations of equation~\!\eqref{quasiPII} as a Hamiltonian system
\begin{align}
H_1^{{\text{qsi-P}}_{\text{I\!I}}} & = \frac{1}{2}\, x^2 - \frac{1}{2}\, y^6 - \frac{1}{4}\, a_3(z)\, y^4 - \frac{1}{3}\, a_2(z)\, y^3 - \frac{1}{2}\, a_1(z)\, y^2 - a_0(z)\, y \,, \label{eq:quasi_H_1} \\[.7ex]
H_2^{{\text{qsi-P}}_{\text{I\!I}}} & = \frac{1}{2}\, x^2 - x\,y^3 - \frac{1}{4}\, b_3(z)\, y^4 - \frac{1}{3}\, b_2(z)\, y^3 - \frac{1}{2}\, b_1(z)\, y^2 - b_0(z)\, y \,, \label{eq:quasi_H_2} \\[.7ex]
H_3^{{\text{qsi-P}}_{\text{I\!I}}} & = \frac{1}{2}\, x^2 - x \left(y^3 + c_3(z) y\right) - \frac{1}{3}\, c_2(z)\, y^3 - \frac{1}{2}\, c_1(z)\, y^2 - c_0(z)\, y \label{eq:quasi_H_3} \,. 
\end{align}

In the following we will follow the same procedure showed in the previous section to compare the systems associated with $H_1^{{\text{qsi-P}}_{\text{I\!I}}}$ and $H_3^{{\text{qsi-P}}_{\text{I\!I}}}$ to evaluate the coefficients $c_i(z)$~\!\eqref{eq:quasi_H_3}. 
\begin{figure}[t]
    \centering
 \includegraphics[width=.999\textwidth]
 {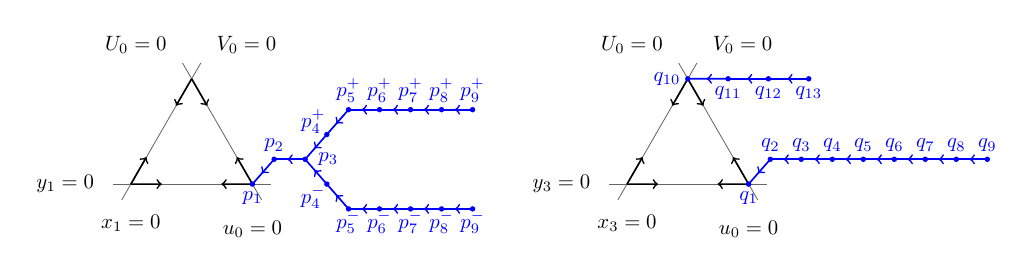}
    \caption{Cascades of blow-ups for the system associated with $H^{{\text{qsi-P}}_{\text{I\!I}}}_1$ on the left and with $H^{{\text{qsi-P}}_{\text{I\!I}}}_3$ on the right.} 
    \label{fig:cascade_quasi_2}
\end{figure}

\begin{figure}[t]
    \centering
   \includegraphics[width=.55\textwidth]{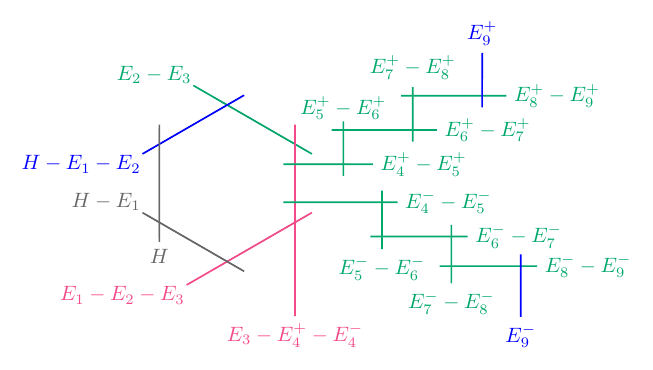}
    \caption{Space of initial conditions for the quasi-Painlev\'e system with $H_1^{{\text{qsi-P}}_{\text{I\!I}}}$.}
    \label{fig:quasi_P_H1}
\end{figure}

%\subsubsection{Blow-up cascades for \texorpdfstring{$H_1^{\,{\text{qsi-P}}_{\text{I\!I}}}$}{} }
The space of initial conditions of the system arising from $H_1^{{\text{qsi-P}}_{\text{I\!I}}}$ was already computed in~\!\cite{KeckerFilipuk}. There the following cascades of blow-ups are given. With the coefficients $a_1(z)$ and $a_3(z)$ fixed as in~\!\eqref{eq:contraints_quasi_PII} we have
\begin{equation} 
\begin{split} \label{eq:blow_up_branch1_quasi_H_1}
    p_1&: (u_0,v_0) = (0\,,0) \,\, \longleftarrow \,\, p_2: (U_1,V_1) = (0\,,0) \,\, \longleftarrow \,\, p_3: (u_2,v_2) = (0\,,0) \\[.7ex]
    & \,\, \longleftarrow \,\, p_4^{+}: (u_3^{+},v_3^{+}) = (0\,,1) \,\, \longleftarrow \,\, p_5^{+}: (u_4^{+},v_4^{+}) = (0\,,0) \,\, \longleftarrow \,\, p_6^{+}: (u_5^{+},v_5^{+}) = \!\Big( 0\,, -\frac{z}{4} \Big) \\[.7ex]
    & \,\, \longleftarrow \,\, p_7^{+}: (u_6^{+},v_6^{+}) = \!\Big( 0\,, -\frac{1}{3}\,a_2(z) \Big) \,\, \longleftarrow \,\, p_8^{+}: (u_8^{+},v_8^{+}) = \!\Big( 0\,, \frac{1}{8} + \frac{z^2}{16} - \frac{c}{2} \Big) \\[.7ex]
    & \,\, \longleftarrow \,\, p_9^{+}: (u_8^{+},v_8^{+}) = \!\Big( 0\,,  \frac{1}{3}\,a_2'(z) +  \frac{z}{4}\, a_2(z) - a_0(z) \Big)\,, 
\end{split} 
\end{equation}
and for the second branch
\begin{equation} 
\begin{split} \label{eq:blow_up_branch2_quasi_H_1}
p_3 &\,\, \longleftarrow\,\, p_4^{-}\colon (u_3^{-},v_3^{-}) = (0\,,-1)\,\, \longleftarrow\,\, p_5^{-}\colon (u_4^{-},v_4^{-}) = (0\,,0)\,\, \longleftarrow p_6^{-}\colon (u_5^{-},v_5^{-}) = \!\Big( 0\,,\frac{z}{4} \Big) \\ 
&\,\, \longleftarrow\,\, p_7^{-}\colon(u_6^{-},v_6^{-}) =  \!\Big( 0\,, \frac{1}{3}\,a_2(z) \Big)\,\, \longleftarrow\,\, p_8^{-}\colon (u_7^{-},v_7^{-}) = \!\Big( 0\,, \frac{1}{8} - \frac{z^2}{16} + \frac{c}{2} \Big) \\ 
&\,\, \longleftarrow\,\, p_9^{-}\colon (u_8^{-},v_8^{-}) = \!\Big( 0\,, \frac{1}{3}\,a_2'(z) -  \frac{z}{4}\, a_2(z) + a_0(z) \Big) \,. 
\end{split}
\end{equation}
The corresponding diagram is depicted in figure~\!\ref{fig:quasi_P_H1}.

With the auxiliary function (obtained in~\!\cite{KeckerFilipuk})
\begin{equation}
W_1^{{\text{qsi-P}}_{\text{I\!I}}}(z) = H_1^{{\text{qsi-P}}_{\text{I\!I}}}(x(z),y(z);z) + \frac{1}{8} \frac{x(z)}{y(z)} + \frac{a_2'(z)}{3} \frac{x(z)}{y(z)^2} + \left( \frac{a_2(z)}{24} + \frac{a_2''(z)}{3} - a_0'(z) \right) \frac{x(z)}{y(z)^4},
\end{equation}
we apply Lemma \ref{log_bounded} to show that the divisors, schematically shown in figure~\!\ref{fig:quasi_P_H1} as green and magenta lines, are inaccessible. Namely, after each blow-up we check that $W_1^{{\text{qsi-P}}_{\text{I\!I}}}$ is infinite at points on the exceptional curve away from any new base points, while the logarithmic derivative of the auxiliary function $W_1^{{\text{qsi-P}}_{\text{I\!I}}}$ is bounded in a neighbourhood of any such point.

\begin{figure}[t]
    \centering
   \includegraphics[width=.5\textwidth]{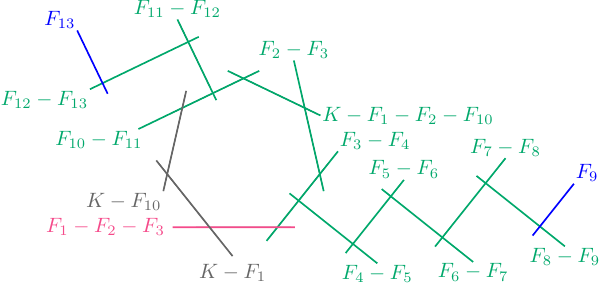}
    \caption{Intersection diagram for the quasi-Painlev\'e system with $H_3^{{\text{qsi-P}}_{\text{I\!I}}}$.}
    \label{fig:quasi_P_H3}
\end{figure}

We now turn to the systems $H_2^{{\text{qsi-P}}_{\text{I\!I}}}$ and $H_3^{{\text{qsi-P}}_{\text{I\!I}}}$, which, although different from the system $H_1^{{\text{qsi-P}}_{\text{I\!I}}}$, result in the same second-order equation~\!\eqref{quasiPII} after elimination of $x(z)$. Both systems $H_2^{{\text{qsi-P}}_{\text{I\!I}}}$ and $H_3^{{\text{qsi-P}}_{\text{I\!I}}}$ have two initial base points in $\mathbb{CP}^2$, located at $q_1\colon (u_0,v_0)=(0\,,0)$ and $q_{10}\colon (U_0,V_0)=(0\,,0)$, with cascades of $9$ blow-ups and $4$ blow-ups required to resolve these points. We here perform the analysis only for the system $H_3^{{\text{qsi-P}}_{\text{I\!I}}}$, as both cases are similar. The base point $q_1$ can be resolved by the following cascade of $9$ blow-ups:
\begin{equation}
\begin{aligned}
    q_1& \colon  (u_0,v_0)=(0\,,0)\,\, \longleftarrow\,\, q_2\colon (U_1,V_1) = (0\,,0)\,\, \longleftarrow\,\, q_3\colon (u_2,v_2) = (0\,,0) \\[2ex]
    & \longleftarrow\,\, q_4\colon (U_3,V_3) = (2,0) \,\,  \longleftarrow\,\, q_5\colon (U_4,V_4)=(0\,,0)\,\, \longleftarrow\,\, q_6\colon (U_5,V_5) = (8\,c_3(z),0) \\[2ex]
    &  \longleftarrow\,\, q_7\colon (U_6,V_6) = \!\Big(\, \frac{8}{3}\,c_2(z) , 0 \Big) %\\    &
   \,\, \longleftarrow\,\, q_8\colon (U_7,V_7) = \!\Big( 8\left(c_1(z) + 8\,c_3^2(z) - 2\,c_3'(z)\right),0 \Big) \\[1.5ex] 
    &   \longleftarrow\,\, q_9\colon (U_8,V_8) = \left( 32 \left(c_0(z) + \frac{4}{3}\,c_2(z)\,c_3(z) -\frac{1}{3}\,c_2'(z) \right) , 0 \right),
\end{aligned}
\end{equation}
the condition for the system after the last blow-up to have algebraic series expansions being 
\begin{equation}
\label{H3QPII_cond1}
    c_3''(z) + \frac{1}{2} c_1'(z) = 0.
\end{equation}
The other base point, $q_{10}$, is resolved via the following cascade of $4$ blow-ups:
\begin{equation}
\begin{split}
q_{10}& \colon  (U_0,V_0) = (0\,,0) \,\, \longleftarrow \,\, q_{11}\colon (U_{10},V_{10}) = \left( -\frac{c_2(z)}{3}, 0 \right)\,\,\longleftarrow \,\, q_{12}\colon (U_{11},V_{11}) = \left( -\frac{c_1(z)}{2} , 0 \right) \\
& \longleftarrow \,\, q_{13}\colon (U_{12},V_{12}) = \left( \frac{1}{3} (c_2(z)\,c_3(z) - c_2'(z) - 3\,c_0(z)) , 0 \right),
\end{split}
\end{equation}
the conditon for the system to have algebraic series expansions here being
\begin{equation}
\label{H3QPII_cond2}
c_1'(z) = 0.
\end{equation}
The two conditions \eqref{H3QPII_cond1} and \eqref{H3QPII_cond2} amount to saying that $c_1$ is a constant and $c_3(z)$ a function at most linear in $z$. In following we will assume that $c_3$ is not constant and normalise it such that $c_3(z)=z$.
The rational surface resulting from the two cascades of blow-ups is schematised in figure~\!\ref{fig:quasi_P_H3}. 

For the Hamiltonian $H_3^{{\text{qsi-P}}_{\text{I\!I}}}$, the auxiliary function needed to show that the divisors of intermediate blow-ups are inaccessible looks as follows:
\begin{equation} \label{eq:W_3_QP_II}
\begin{aligned}
W_3^{{\text{qsi-P}}_{\text{I\!I}}}(z) =&~ H_3^{{\text{qsi-P}}_{\text{I\!I}}}\big(x(z),y(z);z\big) + \frac{1}{2} \frac{x(z)}{y(z)} + \frac{c_2'(z)}{6} \frac{x(z)}{y(z)^2} + \frac{c_2(z)\,c_2'(z)}{9}\,\frac{y(z)}{x(z)} 
\end{aligned}
\end{equation}
After %each of the blow-up process, 
each blow-ups cascade,  
one has to check that $W_3^{{\text{qsi-P}}_{\text{I\!I}}}(z)$ is infinite on the exceptional curve, while ${(W_3^{{\text{qsi-P}}_{\text{I\!I}}})'}/{W_3^{{\text{qsi-P}}_{\text{I\!I}}}}$ is finite in a neighbourhood of any point on the exceptional curve other than any newly arising base points. We did so using Mathematica. In appendix \ref{app:auxiliary_function} we explain how to obtain the terms of the auxiliary function. If the structure of the additional terms is somewhat known, one can also obtain their coefficients by first leaving them as arbitrary values, computing the logarithmic derivative of $W$ and adjusting them for certain terms in this expression to vanish.

\subsection{Comparison of two Hamiltonian systems for \texorpdfstring{${\text{qsi-P}}_{\text{II}}$}{qpII}}
We are now in the position to compare the blow-up diagrams of the Hamiltonians $H_1^{{\text{qsi-P}}_{\text{I\!I}}}$ and $H_3^{{\text{qsi-P}}_{\text{I\!I}}}$ in order to obtain the bi-rational transformation that connects these two systems.

In figure~\ref{fig:quasi_P_H1}, we see that there are $11$ curves of self-intersection number $-2$- and two $-3$-curves. However, there is also a $-1$-curve, denoted $H-E_1-E_2$ in the diagram. This curve is inaccessible by the argument given above and therefore can be blown down. This results in the $-2$-curve $E_2-E_3$ to become a $-1$-curve, $H-E_1-E_3$, which likewise can be blown down. This in turn results in the $-3$-curve $E_3-E_4^{+}-E_4^{-}$ becoming a $-2$-curve, $H-E_1-E_4^{+}-E_4^{-}$. The resulting minimal diagram thus has $11$ $-2$-curves and one $-3$-curve. One can represent this in more abstract form by the intersection diagram below, where nodes represent irreducible components %of the infinity set.
{of the minimal configuration of exceptional curves to be removed from the extended phase space.  }
\begin{equation*}
    \begin{tikzpicture}[scale=.75,baseline={([yshift=-.5ex]current bounding box.center)},scale=0.8]
        \coordinate (o1) at (0,-1);
        \draw (o1) -- coordinate (M) ++(10,0);
        \draw (M) -- ++(0,1) coordinate (N); 
        \foreach \i in {1,...,11} {
        \path[draw=black,fill=white] (\i-1,-1) circle (2.5pt) node[below=.1] {$\delta_{\i}$};
        } 
        \path[draw=black,fill=black] (N) circle (2.5pt) node[above=.1] {$\delta_{12}$};
    \end{tikzpicture} \hspace{5ex}  
    \begin{tabular}{ l l}
    $\delta_1 \colon \mathcal{E}_8^- - \mathcal{E}_9^-$ \quad  &$\delta_7 \colon \mathcal{E}_4^+ - \mathcal{E}_5^+$ \\[.2ex]
    $\delta_2 \colon \mathcal{E}_7^- - \mathcal{E}_8^-$   &$\delta_8 \colon \mathcal{E}_5^+ - \mathcal{E}_6^+$\\[.2ex]
    $\delta_3 \colon \mathcal{E}_6^- - \mathcal{E}_7^-$ &$\delta_9 \colon \mathcal{E}_6^+ - \mathcal{E}_7^+$\\[.2ex]
    $\delta_4 \colon \mathcal{E}_5^- - \mathcal{E}_6^-$   &$\delta_{10} \colon \mathcal{E}_7^+ - \mathcal{E}_8^+$\\[.2ex]
    $\delta_5 \colon \mathcal{E}_4^- - \mathcal{E}_5^-$ &$\delta_{11} \colon \mathcal{E}_8^+ - \mathcal{E}_9^+$\\[.2ex]
    $\delta_6 \colon \mathcal{H}-\mathcal{E}_1-\mathcal{E}_4^+-\mathcal{E}_4^-$ \quad   &$\delta_{12} \colon \mathcal{E}_1 - \mathcal{E}_2 - \mathcal{E}_3$\\
    \end{tabular}
\end{equation*}
%Single circled nodes represent $(-2)$-curves, whereas multiply circled nodes stand for higher self-intersection numbers. 
White nodes represent $-2$-curves, whereas black nodes stand for curves with self-intersection number less than $-2$.
Only connected nodes have non-zero intersection form, indicated by the number of lines (here all equal to $1$). Note that, although the intersection diagram may look a bit like a Dynkin diagram (which it is in the case of the Painlev\'e equations), this does not correspond to any kind of Lie algebra.

First we note that the diagram for $H_1^{{\text{qsi-P}}_{\text{I\!I}}}$ (figure~\!\ref{fig:quasi_P_H1}) an additional $-1$-curve $H-E_1-E_2$ and an additional $-2$-curve $E_2-E_3$. The $-1$-curve can be blown down, which turns $E_1-E_2$ into a new $(-1)$ curve $H-E_1-E_3$, which can be blown down in turn to render the $-3$-curve $E_3-E_4^{+}-E_4^{-}$ a $-2$-curve $H-E_1-E_4^{+}-E_4^{-}$. In this way the intersection diagram becomes minimal and can be matched with the diagram for $H_3^{{\text{qsi-P}}_{\text{I\!I}}}$ (figure~\!\ref{fig:quasi_P_H3}), which is already minimal. However, in order to obtain the change of variables $(x_1,y_1) \mapsto (x_3,y_3)$ we can compare the two non minimal diagrams given by the one represented in figure~\!\ref{fig:quasi_P_H1} and that obtained by blowing up twice the line $F_3-F_4$ in figure~\!\ref{fig:quasi_P_H3}. This is realised by blowing up the $-2$-curve $F_3-F_4$ at an arbitrary point producing the $-3$-curve $F_3-F_4-F_{14}$ intersecting the $-1$-curve $F_{14}$. The latter is blown up as well at an arbitrary point obtaining the $-2$-curve $F_{14}-F_{15}$ intersecting the $-1$-curve $F_{15}$. 

We can compare the Picard lattices associated with the two Hamiltonians 

\vspace*{-1ex}

\begin{equation*}
    \text{Pic}_1^{{\text{qsi-P}}_{\text{I\!I}}}(X) = \text{Span}_{\mathbb{Z}}\,\{\mathcal{H},\mathcal{E}_1,\dots,\mathcal{E}_9^{\pm}\} \simeq \text{Pic}_3^{{\text{qsi-P}}_{\text{I\!I}}}(X) = \text{Span}_{\mathbb{Z}}\,\{\mathcal{K},\mathcal{F}_1,\dots,\mathcal{F}_{15}\}
\end{equation*}
\begin{equation}
\begin{aligned}
    \mathcal{F}_1-\mathcal{F}_2-\mathcal{F}_3&=\mathcal{E}_1-\mathcal{E}_2-\mathcal{E}_3, & \qquad   \mathcal{F}_2-\mathcal{F}_3&=\mathcal{E}_4^--\mathcal{E}_5^-, \\
    \mathcal{F}_3-\mathcal{F}_4-\mathcal{F}_{14}&=\mathcal{E}_3-\mathcal{E}_4^+-\mathcal{E}_4^-, & \qquad  \mathcal{K}-\mathcal{F}_1-\mathcal{F}_2-\mathcal{F}_{10}&=\mathcal{E}_5^--\mathcal{E}_6^-, \\
    \mathcal{F}_4-\mathcal{F}_5&=\mathcal{E}_4^+-\mathcal{E}_5^+, &  \mathcal{F}_{10}-\mathcal{F}_{11} &= \mathcal{E}_6^--\mathcal{E}_7^-, \\
    \mathcal{F}_5-\mathcal{F}_6&=\mathcal{E}_5^+-\mathcal{E}_6^+, &  \mathcal{F}_{11}-\mathcal{F}_{12} &= \mathcal{E}_7^--\mathcal{E}_8^-, \\
    \mathcal{F}_6-\mathcal{F}_7&=\mathcal{E}_6^+-\mathcal{E}_7^+, &  \mathcal{F}_{12}-\mathcal{F}_{13} &= \mathcal{E}_8^--\mathcal{E}_9^-, \\
    \mathcal{F}_7-\mathcal{F}_{8}&=\mathcal{E}_7^+-\mathcal{E}_8^+, &  \mathcal{F}_{14}-\mathcal{F}_{15} &= \mathcal{E}_2-\mathcal{E}_3, \\
    \mathcal{F}_8-\mathcal{F}_{9}&=\mathcal{E}_8^+-\mathcal{E}_9^+, &  \mathcal{F}_{15} &=\mathcal{H}-\mathcal{E}_1-\mathcal{E}_2, 
\end{aligned}
\end{equation}
which allow us to identify the relations between classes of divisors 
\begin{equation}
    \begin{aligned}
        \mathcal{F}_1 &= 2\,\mathcal{H} -\mathcal{E}_1-\mathcal{E}_2-\mathcal{E}_3-\mathcal{E}_4^--\mathcal{E}_5^-,   &\qquad      \mathcal{F}_9 &= \mathcal{E}_9^+,~~\mathcal{F}_{10} = \mathcal{E}_6^-,~~\mathcal{F}_{11} = \mathcal{E}_7^-, \\
        \mathcal{F}_2 &= \mathcal{H} -\mathcal{E}_1-\mathcal{E}_5^-,   &\qquad      \mathcal{F}_{12} &= \mathcal{E}_8^-,~~\mathcal{F}_{13} = \mathcal{E}_9^-,\\
        \mathcal{F}_3 &= \mathcal{H} -\mathcal{E}_1-\mathcal{E}_4^-,  &\qquad \mathcal{F}_{14} &= \mathcal{H} -\mathcal{E}_1-\mathcal{E}_3,      \\
        \mathcal{F}_4 &= \mathcal{E}_4^+,~~\mathcal{F}_5 = \mathcal{E}_5^+,~~ \mathcal{F}_6 = \mathcal{E}_6^+,  &\qquad      \mathcal{F}_{15} &= \mathcal{H} -\mathcal{E}_1-\mathcal{E}_2,\\
       \mathcal{F}_7 &= \mathcal{E}_7^+,~~\mathcal{F}_8 = \mathcal{E}_8^+,     &\qquad     \mathcal{K} &= 3\,\mathcal{H} -2\,\mathcal{E}_1-\mathcal{E}_2-\mathcal{E}_3-\mathcal{E}_4^--\mathcal{E}_5^-.\\
    \end{aligned}
\end{equation}
The map $(x_1,y_1) \mapsto (x_3,y_3)$ is constructed by considering the expressions for the axes $x_3$ and $y_3$ in $H_3^{{\text{qsi-P}}_{\text{I\!I}}}$
\begin{align}
    \{y_3=0\}&\colon \quad \mathcal{K}-\mathcal{F}_1 
    = \mathcal{H} - \mathcal{E}_1 \,, \label{eq:y_3_quasi_H1_H3} \\[.7ex]
    \{x_3=0\}&\colon \quad \mathcal{K}-\mathcal{F}_{10} = 3\,\mathcal{H} -2\,\mathcal{E}_1-\mathcal{E}_2-\mathcal{E}_3-\mathcal{E}_4^--\mathcal{E}_5^- - \mathcal{E}_6^- \,. \label{eq:x_3_quasi_H1_H3} 
\end{align}
From~\!\eqref{eq:y_3_quasi_H1_H3} we recognise a linear expression describing $y_3$ in terms of the coordinates $(x_1,y_1)$ as 
\begin{equation} \label{eq:change_qsi-PII_y3}
    y_3 = A\,x_1 + B\,y_1 + C \,.
\end{equation}
Imposing that the line goes through the point $p_1 \colon (u_0,v_0)=(0\,,0)$, this leads to fix the coefficient $A=0$. The coefficients $B, C$ will be fixed later on. 
From~\!\eqref{eq:x_3_quasi_H1_H3} we have a cubic expression for $y_3$ in $(x_1,y_1)$ as 
\begin{equation} \label{eq:x_3_general_cubic}
    x_3 = a_{30}\,x_1^3 + a_{21}\,x_1^2\,y_1 + a_{12}\,x_1\,y_1^2 + a_{03}\,y_1^3 + a_{20}\,x_1^2 + a_{11}\,x_1\,y_1 + a_{02}\,y_1^2 + a_{10 }\,x_1 + a_{01}\,y_1 + a_{00} \,,
\end{equation}
going through the point $p_1$ with multiplicity $2$ and with simple multiplicity through the points $p_2$, $p_3$, $p_4^-$, $p_5^-$, $p_6^-$, whose coordinates are given in~\!\eqref{eq:blow_up_branch1_quasi_H_1} and~\!\eqref{eq:blow_up_branch2_quasi_H_1}. After imposing these constraints we obtain the resulting expression in terms of the arbitrary coefficients $a_{03}$ and $a_{00}$
\begin{equation} \label{eq:change_qsi-PII_x3}
    x_3 = a_{03} \left(y_1^3 + x_1 + \frac{z}{4}\,y_1\right) + a_{00} \,. 
\end{equation}
We will completely fix the remaining coefficients $B,C,a_{30}$ and $a_{00}$ by requiring that the system given in terms of $x_3(x_1,y_1)$ and $y_3(x_1,y_1)$ is still related to the Hamiltonian $H_3^{{\text{qsi-P}}_{\text{I\!I}}}$ in~\!\eqref{eq:quasi_H_3}
%. In particular, we consider $H_3^{{\text{qsi-P}}_{\text{I\!I}}}$ in~\!\eqref{eq:quasi_H_3} with the coefficients $c_i(z)$ expressed in terms of $a_i(z)$ appearing in~\!\eqref{quasiPII} with constraints~\!\eqref{eq:contraints_quasi_PII}
\begin{equation}
    \begin{cases}
        y_3' = - y_3^3 + x_3  - c_3(z)\, y_3 \\[.7ex]
        x_3' =  x_3\left( 3\,y_3^2+c_3(z) \right) + a_2(z)\,y_3^2 + c_1(z)\, y_3 + c_0(z)  
    \end{cases} \,,
\end{equation}
making use of~\!\eqref{eq:change_qsi-PII_y3},~\!\eqref{eq:change_qsi-PII_x3}, and of the system associated with the Hamiltonian $H_1^{{\text{qsi-P}}_{\text{I\!I}}}$ in~\!\eqref{eq:quasi_H_1}
\begin{equation}
    \begin{cases}
        y_1' =  x_1\\[.7ex]
        x_1' =  3 \, y_1^5 + z \, y_1^3 + a_2(z) \, y_1^2 + \left( \dfrac{z^2}{16} + c \right) y_1 +a_0(z)
    \end{cases} \,.
\end{equation}
This approach will also allow us to relate the coefficients appearing in $H_3^{{\text{qsi-P}}_{\text{I\!I}}}$ to those in $H_1^{{\text{qsi-P}}_{\text{I\!I}}}$. 
Comparing the coefficients of same degree in the variables $x_1$ and $y_1$ we fix the parameters as 
\begin{equation}
  \begin{cases} 
  B = a_{03} = \pm 1 \\  
  C=a_{00} =0 
  \end{cases} \quad \implies \quad      \begin{cases}
        x_3 = \pm \left(y_1^3 + x_1 + \dfrac{z}{4}\,y_1 \right)\\[.7ex] 
        y_3 = \pm \, y_1 
       \end{cases} \,.
\end{equation}
In addition, the set of coefficients in $H_3^{{\text{qsi-P}}_{\text{I\!I}}}$ expressed in terms of those appearing in $H_1^{{\text{qsi-P}}_{\text{I\!I}}}$ with~\!\eqref{eq:contraints_quasi_PII} is 
\begin{equation}
    c_0(z) = \pm \, a_0(z)\,, \quad c_1(z) = \frac{1}{4}+c \,, \quad c_2(z) = \pm \, a_2(z)\,, \quad c_3(z) = \frac{z}{4} \,,
\end{equation}
producing the equations of second order in $y_3$ 
\begin{equation}
    y_3'' = 3 \, y_3^5 + z \, y_3^3 \pm a_2(z) \, y_3^2 + \left( \dfrac{z^2}{16} + c \right) y_3 \pm a_0(z) \,, 
\end{equation}
which can be mapped one onto the other by considering the further transformation $y_3 \to - y_3$. 

Analogously to what we have already observed in the previous section, the form of the map $(x_1,y_1)\mapsto (x_3,y_3)$ depends on the specific choice that we tacitly made between the two branches of the diagram labelling the divisors associated with the two systems, analogously to what happens for the case of $\text{P}_{\text{I\!I}}$. In particular, we could have chosen to identify $F_{13}=E_9^+$ and $F_{9}=E_9^-$, and the curve going through the points $p_4^+, p_5^+,p_6^+$ instead of $p_4^-, p_5^-,p_6^-$. With this choice the change of variables is 
\begin{equation}
    \begin{cases}
        x_3 = \tilde{a}_{03} \left(y_1^3 - x_1 + \dfrac{z}{4}\,y_1 \right) + \tilde{a}_{00}\\[.7ex] 
        y_3 = \tilde{B}\, y_1 + \tilde{C} 
    \end{cases} \,. 
\end{equation}
With the requirement that the system is Hamiltonian with respect to $H_3^{{\text{qsi-P}}_{\text{I\!I}}}$ we have 
\begin{equation}
     \begin{cases} 
  \tilde{B} = \tilde{a}_{03} = \pm i  \\  
  \tilde{C}=\tilde{a}_{00} =0 
  \end{cases} \quad \implies \quad   \begin{cases}
        x_3 = \mp i \left(y_1^3 - x_1 + \dfrac{z}{4}\,y_1 \right)\\[.7ex] 
        y_3 = \pm i\, y_1 
    \end{cases} \,. 
\end{equation}
The coefficients $c_i(z)$ are then determined in terms of $a_i(z)$, as 
\begin{equation}
    c_0(z) = \pm i\, a_0(z)\,, \quad c_1(z) = -\frac{1}{4}+c \,, \quad c_2(z) = \mp i\, a_2(z)\,, \quad c_3(z) = -\frac{z}{4} \,,
\end{equation}
producing the second order ODE in $y_3$ 
\begin{equation}
    y_3'' = 3 \, y_3^5 - z \, y_3^3 \mp i\, a_2(z) \, y_3^2 + \left( \dfrac{z^2}{16} + c \right) y_3 \pm i\, a_0(z) \,. 
\end{equation}
The complete matching with the ODE~\!\eqref{quasiPII} is achieved either considering an inversion of the independent variable $t = -z$ or by mapping the variable $y_3 \to \pm i\, y_3$.

\subsection{Hamiltonian systems of quasi-Painlev\'e-IV type}
\label{quasi-PainleveIV}
The fourth Painlev\'e equation,
\begin{equation}
    y'' = \frac{(y')^2}{2y} + \frac{3}{2} \, y^3 + 4\,z\, y^2 + 2(z^2-\alpha) y + \frac{\beta}{y},
\end{equation}
has been the subject of study by many authors. Several Hamiltonian forms, by Okamoto {\cite{okamoto3}}, Its-Prokhorov~\cite{ItsProkhorov}, Kecker~\!\cite{Kecker2015} and Kajiwara-Noumi-Yamada {\cite{Kajiwara2017}, are known leading to this equation. In~\!\cite{Dzhamay2109.06428}, these Hamiltonians are compared using the geometric approach and relations in the form of bi-rational transformations between them established.
\begin{figure}[t]
    \centering
 \includegraphics[width=.99\textwidth]{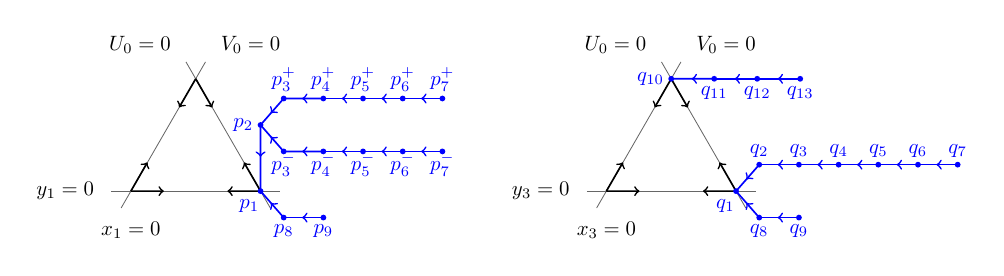}
    \caption{Cascades of blow-ups for the system associated with $H^{{\text{qsi-P}}_{\text{IV}}}_1$ on the left and with $H^{{\text{qsi-P}}_{\text{IV}}}_2$ on the right.} 
    \label{fig:cascade_quasi_4}
\end{figure}

In this section we consider the following equation
\begin{equation}
\label{QPIV}
    y'' = \frac{(y')^2}{2y} + \frac{5}{2} \, y^5 + \alpha_4(z)\,y^4 + \alpha_3(z)\, y^3 + \alpha_2(z)\, y^2 + \alpha_1(z)\, y - \frac{\alpha_0^2(z)}{2y},
\end{equation}
which we propose to be a quasi-Painlev\'e analogue of the fourth Painlev\'e equation.
Like the Painlev\'e-IV equation, a leading order analysis shows that the equation~\!\eqref{QPIV} can have two essentially different types of singularities, namely algebraic poles of the form 
\begin{equation}
\label{leading_QPIV}
    y(z) \sim \frac{\sigma}{\sqrt{2} (z-z_\ast)^{1/2}}\,, \quad \sigma \in \{1,i\}.
\end{equation}
Note that there are four possible values $\sigma \in \{\pm 1,\pm i\}$, but the sign can be absorbed into the choice of branch for $(z-z_\ast)^{-1/2}$.
This leading order behaviour can also be seen when computing the space of initial conditions for this equation, written as a Hamiltonian system. Below we give two different Hamiltonian systems, whose cascades are shown in figure \ref{fig:cascade_quasi_4}, that lead to this equation and compare them using the geometric approach to find a bi-rational change of variables by identifying the rational surfaces of their respective spaces of initial conditions. 
Analogously to the quasi-Painlev\'e II case, we will evaluate the coefficients $\beta_i(z)$ appearing in the Hamiltonian $H^{{\text{qsi-P}}_{\text{IV}}}_2$ \eqref{eq:H_2_qsi-P_IV} in terms of the coefficients $\alpha_i(z)$ in $H^{{\text{qsi-P}}_{\text{IV}}}_1$~\eqref{eq:H_2_qsi-P_IV} after imposing the quasi-Painlev\'e property for the solutions to the related system. 

The first Hamiltonian we consider is
\begin{equation}\label{eq:H_1_qsi-P_IV}
    H^{\,{\text{qsi-P}}_{\text{IV}}}_1(x_1,y_1;z) = \frac{1}{2}\, x_1^2\, y_1 - \frac{1}{2}\, y_1^5 - \frac{\alpha_4(z)}{4}\, y_1^4 - \frac{\alpha_3(z)}{3} \,y_1^3 - \frac{\alpha_2(z)}{2}\, y_1^2 - \alpha_1(z)\, y_1 - \alpha_0(z)\, x_1\,,
\end{equation}
where any terms $x_1^2$ or $x_1 y_1$ could be absorbed by a shift in $x_1$ and $y_1$. The associated system is
\begin{equation}
\label{QPIVsystem1}
\begin{cases}
y_1' = x_1 \,y_1 - \alpha_0(z) \,,\\[.7ex] 
x_1' = \dfrac{5}{2}\,y_1^4 - \dfrac{1}{2}\,x_1^2 + \alpha_4(z) \,y_1^3 + \alpha_3(z) \,y_1^2 + \alpha_2(z) \,y_1 + \alpha_1(z)\,.
\end{cases}    
\end{equation}
Under elimination of $x_1$, equation~\!\eqref{QPIV} is recovered. However, we see that the solution of system~\!\eqref{QPIVsystem1} in fact can exhibit three different types of singular behaviour, with leading orders
\begin{equation}
\label{leading_QPIV_H1_y}
\begin{cases}
    x_1(z) \sim -\dfrac{1}{2(z-z_\ast)},  \\[2ex]  y_1(z) \sim \dfrac{\sigma}{\sqrt{2}(z-z_\ast)^{1/2}}\,, \quad \sigma \in \{1,i\}
\end{cases}\,, \qquad 
\begin{cases}
    x_1(z) \sim \dfrac{2}{z-z_\ast}\,, \\[2ex] y_1(z) \sim \alpha_0(z_\ast)(z-z_\ast)\,,
\end{cases} \,. 
\end{equation}
The first in \eqref{leading_QPIV_H1_y} corresponds to the leading order behaviour~\!\eqref{leading_QPIV} for the branches $\{p_1, \dots, p_7^+\}$ and $\{p_1, \dots, p_7^-\}$ in figure \ref{fig:cascade_quasi_4}. The second term in \eqref{leading_QPIV_H1_y} is the leading-order behaviour for the branch $\{p_1, p_8, p_9\}$ in figure \ref{fig:cascade_quasi_4}, where $x_1$ has a pole and $y_1$ has a zero.

%\subsubsection{Blow-up cascades for \texorpdfstring{$H^{\,{\text{qsi-P}}_{\text{IV}}}_1$}{} }
The system has one base point at the point $p_1\colon (u_0,v_0)=(0\,,0)$ and the regularisation process gives rise to three branches, as in the figure~\!\ref{fig:cascade_quasi_4}. The first branch is 
\begin{equation}
\begin{aligned}
&p_1\,\, \longleftarrow\,\, p_2\colon (U_1,V_1) = (0\,,0)\,\, \longleftarrow\,\, p_3^{+}\colon (u_2^{+},v_2^{+}) = (0\,,1)  \,\, \longleftarrow\,\,  p_4^{+}\colon (u_3^{+},v_3^{+}) = \left( 0\,, -\frac{\alpha_4}{4} \right)\,\,\\[.7ex]
& \longleftarrow\,\, p_5^{+}\colon (u_4^{+},v_4^{+}) = \Big( 0\,,  \frac{3}{2}\,\frac{\alpha_4^2}{16} - \frac{1}{3}\,\alpha_3 \Big) \,\, \longleftarrow\,\, p_6^{+}\colon (u_5^{+},v_5^{+}) = \left( 0\,, f_6(z) \right) \\[2ex]
& \longleftarrow\,\, p_7^{+}\colon (u_6^{+},v_6^{+}) = \left( 0\,, f_7(z) \right)\!,
\end{aligned}
\end{equation}
with 
\begin{align}
    f_6(z) &= \frac{1}{2}\,\frac{\alpha_4'}{4} - \frac{1}{2}\,\alpha_2 + \frac{\alpha_4}{4}\,\alpha_3 -\frac{5}{2}\,\frac{\alpha_4^3}{64} - \alpha_0 \,, \\[.7ex]
    f_7(z) &= 2 \,\frac{\alpha_4}{4}\, \alpha_0 - \alpha_1 +\frac{3}{2}\, \frac{\alpha_4}{4}\, \alpha_2+\frac{1}{6}\, \alpha_3^2 -\frac{5}{2}\, \frac{\alpha_4^2}{16} \, \alpha_3 + \frac{35}{8}\, \frac{\alpha_4^4}{256} + \frac{1}{3}\, \alpha_3' - \frac{5}{2}\, \frac{\alpha_4}{4} \,\frac{\alpha_4'}{4} \,. 
\end{align}
After the $7$th blow-up, we find the following condition to obtain a system with algebraic singularities
\begin{equation}
\label{H7cond}
    \alpha_4''(z) = \frac{3}{16} \,\alpha_4^2(z) \, \alpha_4'(z)- \frac{2}{3} \,\alpha_3(z)\, \alpha_4'(z) - \frac{2}{3}\, \alpha_4(z) \, \alpha_3'(z) + 4\, \alpha_2'(z).
\end{equation}
The second cascade of blow-ups branches off after the second blow-up, where there are two base points $p_3^{\pm}\colon (u_2,v_2) = (0\,,\pm1)$.
\begin{equation}
\begin{split} 
&p_2 \,\, \longleftarrow\,\, p_3^{-}\colon (u_2^{-},v_2^{-}) = (0\,,-1)\,\, \longleftarrow\,\, p_4^{-}\colon (u_3^{-},v_3^{-}) =  \left( 0\,, \frac{\alpha_4}{4} \right) \\[.7ex] 
& \longleftarrow\,\, p_5^{-}\colon (u_4^{-},v_4^{-}) = \Big( 0\,,  -\frac{3}{2}\,\frac{\alpha_4^2}{16} + \frac{1}{3}\,\alpha_3 \Big) \,\,\longleftarrow\,\, p_6^{-}\colon (u_5^{-},v_5^{-}) = \left( 0\,, g_6(z) \right) \\[1ex]  
&\longleftarrow\,\, p_7^{-}\colon (u_6^{-},v_6^{-}) = \left( 0\,, g_7(z)\right)\!,
\end{split}
\end{equation}
with 
\begin{align}
    g_6(z) &=  \frac{1}{2}\,\frac{\alpha_4'}{4} + \frac{1}{2}\,\alpha_2 - \frac{\alpha_4}{4}\,\alpha_3 +\frac{5}{2}\,\frac{\alpha_4^3}{64} - \alpha_0  \,,\\[.7ex]
    g_7(z) &= 2 \,\frac{\alpha_4}{4}\, \alpha_0 + \alpha_1 -\frac{3}{2}\, \frac{\alpha_4}{4}\, \alpha_2-\frac{1}{6}\, \alpha_3^2 +\frac{5}{2}\, \frac{\alpha_4^2}{16} \, \alpha_3 - \frac{35}{8}\, \frac{\alpha_4^4}{256} + \frac{1}{3}\, \alpha_3' - \frac{5}{2}\, \frac{\alpha_4}{4} \,\frac{\alpha_4'}{4} \,. 
\end{align}
Here, after the blow-up of $p_7^{-}$ we find a similar condition as~\!\eqref{H7cond}, but with opposite sign on the right-hand side. Therefore, these two condition decouple into
\begin{equation}
\begin{aligned}
& \alpha_4'' = 0, \\[.7ex] 
&\alpha_4^2 \, \alpha_4 '- \frac{32}{9}\, \alpha_3 \, \alpha_4'- \frac{32}{9}\, \alpha_4\,\alpha_3'+\frac{64}{3}\, \, \alpha_2' = \frac{d}{dz} \left( 3\,\alpha_4^3 - \frac{32}{9}\, \alpha_4\,\alpha_3 + \frac{64}{3}\,\alpha_2 \right) = 0 \,. 
\end{aligned} 
\end{equation}
This means that $\alpha_4$ is at most linear, $\alpha_4(z) = a\,z +b$, and in case $a\neq 0$ we shift and re-scale the expression to obtain 
\begin{equation} \label{eq:quasiP_property_qsi-PIV}
    \alpha_4(z) = 4\,z\,, \quad \alpha_2(z) = c + \frac{2\,z}{3} \, \alpha_3(z) - z^3\,,
\end{equation}
where $c \in \mathbb{C}$ is a constant.  
Lastly, there is a third cascade of blow-ups, given below:
\begin{equation}
p_1 \quad \longleftarrow \quad p_8\colon (u_1,v_1) = (0\,,0) \quad \longleftarrow \quad p_9\colon (u_8,v_8) = (0\,,2\,\alpha_0),
\end{equation}
which gives the additional condition that $\alpha_0$ is a constant. Note that, after the blow-up of $p_9$ we find a regular system, meaning that the singularities arising on the exceptional curve $E_9$ are ordinary poles, whereas singularities arising on either of the exceptional curves $E_7^{\pm}$ are algebraic poles of square-root type. Thus we have an example of a system with mixed types of poles, ordinary as well as algebraic ones. { While an equation of algebro-Painlev\'e type (a stronger requirement than the quasi-Painlev\'e property) with mixed types of poles was already considered in~\cite{FilipukStokes_2023}, this is obtained somewhat artificially, by a simple algebraic change of variables $w=y^2$, from the fourth Painlev\'e equation. Such transformations can always be applied to Painlev\'e equations and lead to equations of algebro-Painlev\'e type. In that way, the system presented here are more general in the way that it is of (genuine) quasi-Painlev\'e type, i.e.\ not transformable to one of the Painlev\'e equations.
}
In system~\eqref{QPIVsystem1}, the coefficients $\alpha_1(z)$ and $\alpha_3(z)$ are still arbitrary analytic functions, so that the general equation~\!\eqref{QPIV} with the quasi-Painlev\'e property is 
\begin{equation}
y'' = \frac{(y')^2}{2\,y} + \frac{5}{2}\, y^5 + 4\,z\,y^4 + \alpha_3(z)\, y^3 + \left(c + \frac{2\,z}{3}\,\alpha_3(z) - z^3\right) y^2 + \alpha_1(z)\, y - \frac{\alpha_0^2}{2\,y},
\end{equation}
where $c \in \mathbb{C}$ and $\alpha_1(z)$, $\alpha_3(z)$ are analytic (but otherwise arbitrary) functions.
\begin{figure}[t]
    \centering
  \includegraphics[width=.4\textwidth]{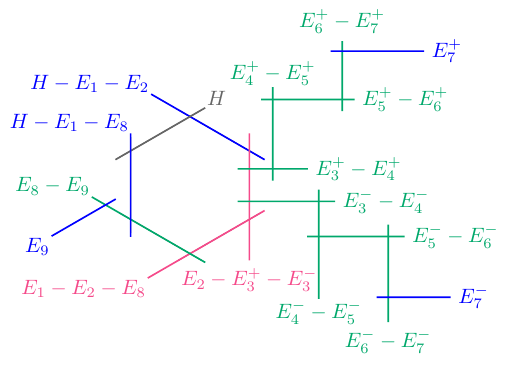}
    \caption{Intersection diagram for the quasi-Painlev\'e system with $H_1^{{\text{qsi-P}}_{\text{IV}}}$.}
    \label{fig:quasi_PIV_H1}
\end{figure}

The auxiliary function built from the Hamiltonian $H_1^{{\text{qsi-P}}_{\text{I\!V}}}$, to show that the exceptional curves of all intermediate blow-ups (magenta and green lines, as well as the blue line $H-E_1-E_2$ in figure~\!\ref{fig:quasi_PIV_H1}) are inaccessible by the solution, is given by 
\begin{equation}
W_1^{{\text{qsi-P}}_{\text{I\!V}}} = H_1^{{\text{qsi-P}}_{\text{I\!V}}} + \frac{1}{2} \frac{x(z)\, y(z)}{y(z) + A_1(z)}\,, \qquad A_1(z)= 3\,z-\frac{2}{3}\,\alpha_3'(z)\,.
\end{equation}
The diagram is depicted in figure \ref{fig:quasi_PIV_H1} and is not minimal. In order to produce the minimal version we could blow down the $-1$-curve $H-E_1-E_2$, so that instead of the $-3$-curve $E_2-E_3^+-E_3^-$ we obtain a $-2$-curve. The resulting diagram is then minimal and we can replicate the procedure followed in the previous section to determine the diagram given in terms of divisors associated with the system.

A second Hamiltonian system giving rise to equation~\!\eqref{QPIV} is given by
\begin{equation}\label{eq:H_2_qsi-P_IV}
    H^{\,{\text{qsi-P}}_{\text{IV}}}_2(x_2,y_2;z) = \frac{1}{2}\,x_2^2\, y_2 -x_2 \left(y_2^3 +\frac{\beta_4(z)}{4}\, y_2^2 +\frac{\beta_3(z)}{3}\, y_2 +\beta_0(z) \right)-\frac{\beta_2(z)}{2}\, y_2^2 - \beta_1 (z) \,y_2 \,,
\end{equation}
\begin{equation}
\begin{cases}
y_2' = x_2\, y_2 - y_2^3 - \dfrac{\beta_4(z)}{4} \,y_2^2 - \dfrac{\beta_3(z)}{3} y_2 - \beta_0(z) \\[2ex] 
x_2' = 3\,x_2\, y_2^2 + \dfrac{\beta_4(z)}{2} \,x_2 \,y_2 + \dfrac{\beta_3(z)}{3} \,x_2 + \beta_2(z) \,y_2 + \beta_1(z) - \dfrac{1}{2} \,x_2^2 
\end{cases}
\end{equation}
with the relation between the coefficients $\beta_i(z)$ and $\alpha_i$ being revealed further below. 
The three types of leading-order behaviour at movable singularities for this system are

\vspace*{-4ex}

\begin{equation} \label{eq:behaviour_qsi-P_IV_H2}
\begin{cases}
x_2(z) \sim \dfrac{-1}{z-z_\ast} \\[2.5ex] y_2(z) \sim \dfrac{i}{\sqrt{2}(z-z_\ast)^{1/2}} 
\end{cases},\quad \begin{cases}
x_2(z) \sim \dfrac{2}{z-z_\ast} \\[2.5ex] y_2(z) \sim \beta_0(z_\ast) (z-z_\ast)\,     
\end{cases},\quad \begin{cases}
x_2(z) \sim -\dfrac{\beta_2(z_\ast)}{\sqrt{2}} (z-z_\ast)^{1/2} \\[2.5ex] y_2(z) \sim \dfrac{1}{\sqrt{2}(z-z_\ast)^{1/2}}\,
\end{cases}. 
\end{equation}
The first term in \eqref{eq:behaviour_qsi-P_IV_H2} refers to the leading order of the behaviour of the solutions in the neighbourhood of the singularity $z_\ast$ for the branch $\{q_1, \dots, q_7\}$ in figure \ref{fig:cascade_quasi_4}. The second and the third term in \eqref{eq:behaviour_qsi-P_IV_H2} account for the branch $\{q_1,q_8,q_9\}$ and $\{q_{10},q_{11},q_{12}\}$ respectively. 
The sequences of blow-ups to perform in order to regularise the system described by $H^{\,{\text{qsi-P}}_{\text{IV}}}_2$ are explicitly given in the following. The system shows two base points, at 

\vspace*{-4ex}

\begin{equation} \label{eq:quasi_P4_H2_base_points}
    q_1\colon (u_0,v_0) = (0\,,0)\,, \qquad q_{10}\colon (U_0,V_0) = (0\,,0) \,,
\end{equation}
from which the cascade of blow-ups are in the following. The first branch is  

\vspace*{-4ex}

\begin{align}
\begin{split}
 q_1& \longleftarrow\,\, q_2\colon (U_1,V_1) = (0\,,0)\,\, \longleftarrow\,\, q_3\colon (U_2,V_2) = (2,0)  \\[.5ex] 
& \hspace*{-2ex}\longleftarrow\,\, q_4\colon (U_3,V_3) = (\beta_4,0) \,\, \longleftarrow\,\, q_5\colon (U_4,V_4) = \left(  \frac{1}{2}\,\beta_4^2 + \frac{8}{3}\,\beta_3 \,,\,0\right)\\[.5ex] 
& \hspace*{-2ex} \longleftarrow\,\, q_6\colon (U_5,V_5) = \left( h_6(z) \,,\,0 \right) \,\, \longleftarrow\,\, q_7\colon (U_6,V_6) = \left( h_7(z)\,,\,0 \right) \,, 
\end{split}
\end{align}
with

\vspace*{-6ex}

\begin{align}
    h_6(z) &= 4 \,\beta_4 \,\beta_3 +\frac{\beta_4^3}{4}+4 \,\beta_2+16\, \beta_0 -2 \,\beta_4' \,,\\[.7ex]
    \begin{split} 
    h_7(z) &=  32\, \beta_0  \, \beta_4+\frac{1}{8}\, \beta_4^4+4\, \beta_4^2 \,\beta_3+\frac{64}{3}\, \beta_3^2+6\, \beta_4 \,\beta_2+16\, \beta_1-3\, \beta_4 \,\beta_4'-\frac{32}{3}\, \beta_3'\,.
    \end{split}
\end{align}
The second branch starting from $q_1$ is 
\vspace*{-2ex}

\begin{align} 
q_1 &\longleftarrow\,\, q_8: (u_1,v_1) = (0\,,0)\,\, \longleftarrow\,\, q_9: (u_8,v_8) = (0\,,2\beta_0) \,,  
\end{align}
and as before, the condition that the system after the blow-up of $q_9$ is regularisable is $\beta_0' = 0$, i.e.\ that $\beta_0$ is a constant. In this case the system is in fact regular, i.e.\ in the coordinates $(u_9,v_9)$ it has analytic solutions on the exceptional curve $F_9$, which translate into poles in the variable $x$, whereas $y$ has a simple zero at such points. The branch originated at the point $q_{10}$ is 
\vspace*{-4ex}

\begin{align}
\begin{split} 
\!\!q_{10} & \longleftarrow\, q_{11}\colon (U_{10},V_{10}) = (0\,,0)\,\, \longleftarrow\,\, q_{12}\colon (U_{11},V_{11}) = \!\Big( -\frac{\beta_2(z)}{2}\, , 0 \Big) \, \\ 
&\longleftarrow\,\, q_{13}\colon (U_{12},V_{12}) = \Big( \frac{\beta_4\, \beta_2}{8}\,  - \beta_1 , 0 \Big) \,.
\end{split}
\end{align}
The condition on the system to be regularisable after the blow-up of $q_{13}$ is $\beta_2'(z) = 0$, i.e.\ $\beta_2$ has to be a constant as well. Together with the above condition this implies that $\beta_4''(z) = 0$, that is $\beta_4(z) = az + b$ is at most linear, and in the case $a \neq 0$ we shift and re-scale $z$ so that $\beta_4(z) = 4z$.

\begin{figure}[t]
    \centering
    \includegraphics[width=.5\textwidth]{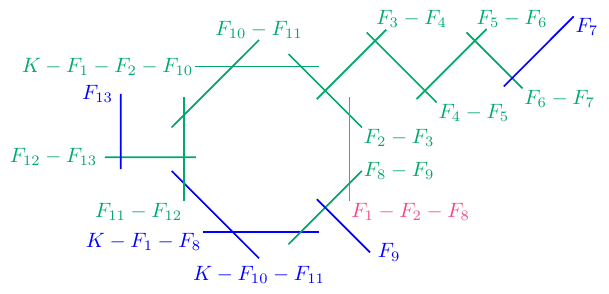}
    \caption{Space of initial conditions for the quasi-Painlev\'e system with $H_2^{{\text{qsi-P}}_{\text{IV}}}$. }
    \label{fig:quasi_PIV_H3}
\end{figure}

With these conditions satisfied, we can derive the second-order equation obtained from the Hamiltonian $H_2^{{\text{qsi-P}}_{\text{IV}}}$,
\vspace*{-2ex}

\begin{equation}
\begin{aligned}
y'' &= \frac{(y')^2}{2 y}+\frac{5}{2}\,y^5+4 \,z\,y^4+\Big(\beta_3  +\frac{3}{2} \,z^2 \Big)y^3
+\Big(2 \beta_0 + \frac{2}{3} \,z\,\beta_3 +\beta_2 -1\Big) y^2 \\
&~~
+\Big( \beta_1 + z \,\beta_0  -\frac{1}{3} \,\beta_3 '+\frac{1}{18} \,\beta_3 ^2 \Big)  y -\frac{\beta_0 ^2}{2 y}\,. 
\end{aligned} 
\end{equation}
The auxiliary function for this system to show that intermediate exceptional curves are inaccessible is 
\begin{equation}
W_2^{{\text{qsi-P}}_{\text{I\!V}}} = H_2^{{\text{qsi-P}}_{\text{I\!V}}} + \frac{1}{2} \frac{x(z)\, y(z)}{y(z) + B_1(z)}, \qquad B_1(z) = z - \frac{2}{3} \,\beta_3'(z)\,.
\end{equation}
As before, the auxiliary function has a bounded logarithmic derivative along all the inaccessible divisors, and it is bounded along the final exceptional curves at the end of each cascade of blow-ups. 

\subsection{Comparison of two Hamiltonian systems for \texorpdfstring{${\text{qsi-P}}_{\text{IV}}$}{qpIV}}
Following the approach illustrated for quasi-Painlev\'e II, we will evaluate the coefficients $\beta_i(z)$ appearing in $H_2^{{\text{qsi-P}}_{\text{IV}}}$ in terms of the coefficients $\alpha_i(z)$ appearing in $H_1^{{\text{qsi-P}}_{\text{IV}}}$. To do so, we construct the Picard lattices associated with the two systems and establish the equivalences between {the irreducible components of the inaccessible divisor}. 
The graph in figure~\!\ref{fig:quasi_PIV_H3} representing the rational surface for $H_2^{{\text{qsi-P}}_{\text{IV}}}$ is minimal.
\begin{equation*}
    \begin{tikzpicture}[scale=.75,baseline={([yshift=-.5ex]current bounding box.center)},scale=0.8]
        \coordinate (o1) at (0,-1);
        \draw (o1) -- ++(4,0) coordinate (M);
        \draw (M) -- ++(4,0);
        \draw (M) -- ++(0,1) coordinate (N) -- ++(0,1) coordinate (P); 
        \foreach \i in {1,...,9} {
        \path[draw=black,fill=white] (\i-1,-1) circle (2.5pt) node[below=.1] (\i){$\delta_{\i}$};
        } 
        \path[draw=black,fill=black] (N) circle (2.5pt) node[right=.1] {$\delta_{10}$};
        \path[draw=black,fill=white] (P) circle (2.5pt) node[above=.1] {$\delta_{11}$};
    \end{tikzpicture} \hspace{5ex}  {
    \begin{tabular}{ l l l}
    $\delta_1 \colon \mathcal{F}_{12} - \mathcal{F}_{13}$   & $\delta_6 \colon \mathcal{F}_3-\mathcal{F}_4$ & $\delta_{11} \colon \mathcal{F}_8-\mathcal{F}_9$ \\[.7ex]
    $\delta_2 \colon \mathcal{F}_{11} - \mathcal{F}_{12}$   & $\delta_7 \colon \mathcal{F}_4 - \mathcal{F}_5$\\[.7ex]
    $\delta_3 \colon \mathcal{F}_{10} - \mathcal{F}_{11}$ & $\delta_8 \colon \mathcal{F}_5 - \mathcal{F}_6$\\[.7ex]
    $\delta_4 \colon \mathcal{K}-\mathcal{F}_1-\mathcal{F}_2 - \mathcal{F}_{8} \quad$ & $\delta_9 \colon \mathcal{F}_6 -\mathcal{F}_7$ \\[.7ex]
    $\delta_5 \colon \mathcal{F}_2 - \mathcal{F}_3$ & $\delta_{10} \colon \mathcal{F}_1-\mathcal{F}_2 -\mathcal{F}_8 $ 
    \end{tabular}}
\end{equation*}
In order to give the change of variables $(x_1,y_1)\mapsto (x_2,y_2)$ we will consider the comparison between non minimal surfaces. In particular, we blow up the line $F_2-F_3$ in figure~\!\ref{fig:quasi_PIV_H3} at an arbitrary point to obtain the $-3$-curve $F_2-F_3-F_{14}$ and we can establish the relations among classes of divisors 
\begin{equation}
\begin{aligned}
    \mathcal{F}_1-\mathcal{F}_2-\mathcal{F}_8&=\mathcal{E}_1-\mathcal{E}_2-\mathcal{E}_8, & \hspace{8ex}   \mathcal{F}_2-\mathcal{F}_3-F_{14}&=\mathcal{E}_2-\mathcal{E}_3^+-\mathcal{E}^-, \\
    \mathcal{F}_3-\mathcal{F}_4&=\mathcal{E}_3^--\mathcal{E}_4^-, & \qquad  \mathcal{K}-\mathcal{F}_1-\mathcal{F}_2-\mathcal{F}_{10}&=\mathcal{E}_3^+-\mathcal{E}_4^+, \\
    \mathcal{F}_4-\mathcal{F}_5&=\mathcal{E}_4^--\mathcal{E}_5^-, &  \mathcal{F}_{10}-\mathcal{F}_{11} &= \mathcal{E}_4^+-\mathcal{E}_5^+, \\
    \mathcal{F}_5-\mathcal{F}_6&=\mathcal{E}_5^--\mathcal{E}_6^-, &  \mathcal{F}_{11}-\mathcal{F}_{12} &= \mathcal{E}_5^+-\mathcal{E}_6^+, \\
    \mathcal{F}_6-\mathcal{F}_7&=\mathcal{E}_6^--\mathcal{E}_7^-, &  \mathcal{F}_{12}-\mathcal{F}_{13} &= \mathcal{E}_6^+-\mathcal{E}_7^+, \\
    \mathcal{F}_8-\mathcal{F}_{9}&=\mathcal{E}_8-\mathcal{E}_9, &  \mathcal{F}_{14} &= \mathcal{H}-\mathcal{E}_1-\mathcal{E}_2, \\ 
\end{aligned}
\end{equation}
yielding the following 
\begin{equation}
    \begin{aligned}
        \mathcal{F}_1 &= \mathcal{H} -\mathcal{E}_2-\mathcal{E}_3^+, \qquad &\qquad \mathcal{F}_9 &= \mathcal{E}_9,~~\mathcal{F}_{10} = \mathcal{E}_4^+,~~\mathcal{F}_{11} = \mathcal{E}_5^+, \\
        \mathcal{F}_2 &= \mathcal{H} -\mathcal{E}_1-\mathcal{E}_3^+,    &\qquad \mathcal{F}_{12} &= \mathcal{E}_6^+,~~\mathcal{F}_{13} = \mathcal{E}_7^+,\\
        \mathcal{F}_3 &= \mathcal{E}_3^-,~~\mathcal{F}_4 = \mathcal{E}_4^- ,~~\mathcal{F}_5 = \mathcal{E}_5^-,  &\hspace{8ex}  \mathcal{F}_{14} &= \mathcal{H} -\mathcal{E}_1-\mathcal{E}_2,\\
        \mathcal{F}_6 &= \mathcal{E}_6^-,~~\mathcal{F}_7 = \mathcal{E}_7^-,~~\mathcal{F}_8 = \mathcal{E}_8,  &\qquad \mathcal{K} &= 2\,\mathcal{H} -\mathcal{E}_1-\mathcal{E}_2-\mathcal{E}_3^+.\\
    \end{aligned}
\end{equation}
As before the map $(x_1,y_1) \mapsto (x_2,y_2)$ is constructed by referring to the axes $x_2$ and $y_2$ in $H_2^{{\text{qsi-P}}_{\text{IV}}}$
\begin{align}
    \{y_2=0\}&\colon \quad \mathcal{K}-\mathcal{F}_1-\mathcal{F}_8
    = \mathcal{H} - \mathcal{E}_1-\mathcal{E}_8 \,, \label{eq:y_2_quasiIV_H1_H2} \\[.7ex]
    \{x_2=0\}&\colon \quad \mathcal{K}-\mathcal{F}_{10}-\mathcal{F}_{11} = 2\,\mathcal{H} -\mathcal{E}_1-\mathcal{E}_2-\mathcal{E}_3^+-\mathcal{E}_4^+-\mathcal{E}_5^+ \,. \label{eq:x_2_quasiIV_H1_H2} 
\end{align}
From~\!\eqref{eq:y_2_quasiIV_H1_H2} we recognise a linear expression describing $y_2$ in terms of the coordinates $(x_1,y_1)$ as 
\begin{equation} \label{eq:change_qsi-PIV_y3}
    y_2 = A\,x_1 + B\,y_1 + C \,,
\end{equation}
and going through the points $p_1 \colon (u_0,v_0)=(0\,,0)$ and $p_8 \colon (u_1,v_1)=(0\,,0)$, yielding 
\begin{equation}
    A=C=0 \,, \quad \implies \quad  y_2 = B\,y_1\,. 
\end{equation}
The relation~\!\eqref{eq:x_2_quasiIV_H1_H2} leads to the transformation
\begin{equation}
    x_2 = a_{20}\,x_1^2 + a_{11}\, x_1\,y_1 + a_{02}\,y_1^2 + a_{10}\,x_1 + a_{01}\,y_1 + a_{00} \,,
\end{equation}
going through the points $p_1, p_2, p_3^+,p_4^+,p_5^+$,  fixing some of the constraints to 
\begin{equation}
    x_2 = a_{02}\left(y_1^2 + x_1 +z\,y_1 +\frac{\alpha_3(z)}{3}-\frac{z^2}{2}\right) ,
\end{equation}
where we made use of~\!\eqref{eq:quasiP_property_qsi-PIV}. To completely determine the coordinate transformation we impose that the coordinates satisfy the Hamiltonian equations, and we get 
\begin{equation}
    \begin{cases}
        a_{02}=1 \\
        B= \pm 1
    \end{cases} \quad \implies \quad 
    \begin{cases}
        x_2 = y_1^2 + x_1 +z\,y_1 +\dfrac{\alpha_3(z)}{3}-\dfrac{z^2}{2} \\
        y_2 = \pm y_1
    \end{cases} \,,
\end{equation}
with the coefficients of $H_2^{{\text{qsi-P}}_{\text{I\!V}}}$ 
\begin{equation}
\begin{aligned}
    \beta_0(z) &= \pm \alpha_0(z)\,, \qquad \beta_1(z) = z \left(1+\alpha_0 + \frac{z^3}{8} \right) + \frac{z^2}{2}\,\frac{\alpha_3(z)}{3}-\frac{1}{2}\left(\frac{\alpha_3(z)}{3} \right)^2 + \frac{\alpha_3(z)'}{3} \,, \\[.7ex] 
    \beta_2(z) &= c \pm 2\, \alpha_0(z) \mp 1 \,, \qquad \beta_3(z) = \alpha_3(z) -\frac{3\,z^2}{2} \,, \qquad \beta_4(z) = \pm 4z \,. 
\end{aligned}
\end{equation}

\section{Conclusion}
We have shown how, using the geometric approach of comparing the { irreducible components of the inaccessible divisor} obtained in the { regularisation } process of Hamiltonian systems { via blow-ups}, one can find explicitly bi-rational coordinate transformations between two such related systems. We have studied various systems for which the leading order behaviour of the solutions, is proportional to a term $(z-z_\ast)^{-1/k}$, where $k \in \{1,2\}$. We have seen examples where different cascades of blow-ups lead to a different leading order behaviour, which is also dictated by the form of the $2$-form in coordinates in the final chart of a blow-up cascade. We have thus encountered equations which have solutions with mixed types of singularities, e.g.\ $(z-z_\ast)^{-1}$ and $(z-z_\ast)^{-1/2}$ for the systems in section~\!\ref{quasi-PainleveIV}.

The blow-up process therefore serves two main purposes. Firstly, starting from a given equation or (polynomial) Hamiltonian with analytic functions as coefficients, it allows us to obtain conditions that are necessary for the equations to have the quasi-Painlev\'e property. These conditions are obtained in the systems after the final blow-up in a cascade, where they allow an additional cancellation in the rational expressions of the vector field so that this can be integrated to give algebraic singularities. These are the same conditions that one obtains from a ``quasi-Painlev\'e'' test: by inserting formal Puiseux series solutions into the equations one can compute a recurrence relation for the coefficients that can be used to determine all coefficients apart from one, where the recurrence breaks down. In this case, the resonance condition that needs to be satisfied to give algebraic series is equivalent to the conditions obtained from one of the blow-up cascades of the equation. Otherwise, one is forced to introduce logarithmic terms $\log(z-z_\ast)$ in the series. 

{
Secondly, by computing the space of initial conditions through blow-ups, one finds a rational surface on which the equations are regularised after a further change of dependent and independent variables. Going further, one could write down a global atlas for the space of initial conditions as was done in~\cite{FilipukStokes} with a $2$-form $\omega = dy \wedge dx = u_i^a \, du_i \wedge dv_i$, for each chart $(u_i,v_i)$ at the end of a blow-up cascade. With the resonance conditions satisfied, these would also render the Hamiltonian with respect to these $2$-forms to be polynomial in the respective coordinates $(u_i,v_i)$ of the chart. In the cases we have considered here, the rational surface, or rather the configuration of the inaccessible divisor one removes from the blown-up space, defines the Hamiltonian structure of these equations, in the way that two systems that are related via a bi-rational transformation have the same minimal blow-up diagram. By comparing the irreducible components of the inaccessible divisor of two systems one can explicitly obtain the change of variables that transforms one system into the other. Lastly, we note that in a forthcoming paper, we classify more systematically Hamiltonian systems of quasi-Painlev\'e type with quartic Hamiltonian, and relate the systems studied in the present paper with certain subcases found in this classification~\cite{MDAKec2}.
}
\vspace*{-2ex}

\section*{Acknowledgements}
We would like to thank Anton Dzhamay, Alexander Stokes and Giorgio Gubbiotti for helpful discussions. The authors of this work were supported through the UKRI / EPSRC New Investigator Award EP/W012251/1. 
We would also like to thank the anonymous referee for a thorough review of the article and many helpful suggestions to improve the overall quality of the paper.

\appendix

\vspace*{-1ex}

\section{Change of variables with the geometric approach}
We give the details of the computation to determine the change of variables to compare the Hamiltonian systems associated with $H_1^{\text{P}_{\text{I\!I}}}$ and $H_4^{\text{P}_{\text{I\!I}}}$. From section~\!\ref{P2review} we recall the equivalence between classes of divisors for the axes
\begin{equation} \label{eq:app_constraints_divisors}
\begin{aligned}
    \mathcal{K}-\mathcal{F}_4 &= \mathcal{H}-\mathcal{E}_1, \\[.7ex]
    \mathcal{K}-\mathcal{F}_1-\mathcal{F}_2 &= 2\,\mathcal{H} - \mathcal{E}_1 - \mathcal{E}_2 - \mathcal{E}_3^{-} - \mathcal{E}_4^{-} - \mathcal{E}_5^{-},
\end{aligned}    
\end{equation}
from which we get the following expressions in the coordinates $(x_1(z),y_1(z))$ as functions of $(x_4(z),y_4(z))$
\vspace*{-8ex}

\begin{align}
   x_1 &\colon \quad   y_1 = A\,x_4 + B\,y_4 + C = 0 \,, \label{eq:app_line} \\[.7ex]
   y_1 &\colon \quad   x_1 = a \,x_4^2 + b \,x_4\,y_4 + c \,y_4^2 + d\,x_4 + e\,y_4 + f = 0\,.  \label{eq:app_conic}
\end{align}
\vspace*{-4ex}

\noindent
To completely determine the transformation of variables $(x_1(z),y_1(z)) \mapsto (x_4(z),y_4(z))$ we must fix the free parameters appearing in the expressions, by imposing the constraints coming from~\!\eqref{eq:app_constraints_divisors}. In particular, the line~\!\eqref{eq:app_line} goes through the point $q_1\colon (u_0,v_0)=(0\,,0)$, and the conic~\!\eqref{eq:app_conic} goes through the points $q_1$, $q_2$, $q_3^{-}$, $q_4^{-}$ and $q_5^{-}$ with coordinates 
\vspace*{-4ex}

\begin{equation*}
    q_2 \colon (U_1,V_1) = (0\,,0)\,, ~~ q_3^- \colon (u_2^-,v_2^-) = \left(0\,,\frac{3}{2}\right)\,, ~~ 
    q_4^- \colon (u_3^-,v_3^-) = \left(0\,,0\right)\,, ~~  q_5^- \colon (u_4^-,v_4^-) = \left(0\,,\frac{3}{8}\,z\right)\,. 
\end{equation*}
\vspace*{-4ex}

\noindent
Using the definition of variables $(u_0,v_0)$ in the system associated with $H_4^{\text{P}_{\text{I\!I}}}$, $\big(x_4 ={1}/{u_0},\,y_4 = {v_0}/{u_0}\big)$\,,
we get the following constraint for~\eqref{eq:app_line} 
\vspace*{-1ex}

\begin{equation}
    \begin{aligned}
        x_1&\colon A +  B\, v_0 + C \, u_0 = 0 \,, \\
        q_1&\colon (u_0,v_0)=(0\,,0) \in { \{ x_1 = 0\}}  \quad \implies \quad  A = 0\,. 
    \end{aligned}
\end{equation}
\vspace*{-2ex}

\noindent
The conic~\eqref{eq:app_conic} goes through the same point 
\vspace*{-1ex}

\begin{equation}
    \begin{aligned}
        y_1&\colon a + b\,v_0 + c\, v_0^2 + d\, u_0 + e\, v_0\, u_0 + f\, u_0^2  = 0 \,, \\
        q_1&\colon (u_0,v_0)=(0\,,0) \in { \{y_1 = 0\}} \quad \implies \quad  a = 0\,. 
    \end{aligned}
\end{equation}
We consider the inverse of the change of coordinates yielding the chart $(U_1,V_1)$, i.e.\ $(u_0 = U_1\,V_1,\,v_0 = V_1)$, obtaining the following expression for the conic 
\vspace*{-1ex}

\begin{equation}
    \begin{aligned}
        y_1&\colon b\,V_1 + c\, V_1^2 + d\, U_1\,V_1 + e\, U_1\, V_1^2 + f\, U_1^2\,V_1^2  = 0 \,,
    \end{aligned}
\end{equation}
from which we extract part of the set of points defining the exceptional curve $E_1\colon \{u_1=0\} \cup \{V_1=0\} $ to get 
\vspace*{-1ex}

\begin{equation}
    \begin{aligned}
        y_1&\colon b + c\, V_1 + d\, U_1 + e\, U_1\, V_1 + f\, U_1^2\,V_1  = 0 \,. 
    \end{aligned}
\end{equation}
We impose the second constraint on the curve 
\vspace*{-1ex}

\begin{equation}
    \begin{aligned}
        y_1&\colon b + c\, V_1 + d\, U_1 + e\, U_1\, V_1 + f\, U_1^2\,V_1  = 0 \,, \\[1.5ex]
        q_2&\colon (U_1,V_1)=(0\,,0) \in { \{y_1 =0\}} \quad \implies \quad  b = 0\,. 
    \end{aligned}
\end{equation}
We iterate the procedure for the remaining constraints. The next change of variables is 
\vspace*{-1ex}

\begin{equation}
    U_1 = u_2^- \equiv u_2 \,, \qquad V_1 = u_2^-\,v_2^- \equiv u_2\,v_2 \,,
\end{equation}
and we can apply the next constraint after having select the exceptional curve $E_2 \colon \{ u_2 = 0  \}\cup \{V_2 = 0\}$
\vspace*{-1ex}

\begin{equation}
    \begin{aligned}
        y_1&\colon c\, v_2 + d + e\, u_2\,v_2 + f\, u_2^2\,v_2  = 0 \,, \\[1.5ex]
        q_3^-&\colon (u_2,v_2)=\left(0\,,\frac{3}{2}\right) \in { \{ y_1 =0\}}  \quad \implies \quad  d = -\frac{3}{2}\,c \,. 
    \end{aligned}
\end{equation}
The new chart is given by 
\vspace*{-2ex}

\begin{equation}
    u_2 = u_3^- \equiv u_3 \,, \qquad v_2 = u_3^-\,v_3^-+\frac{3}{2} \equiv u_3\,v_3+\frac{3}{2}\,,
\end{equation}
from which we get the following 
\vspace*{-1ex}

\begin{equation}
    \begin{aligned}
        y_1&\colon  \frac{3}{2}\, e+ \frac{3}{2}\, f\, u_3 +c\,v_3 + e\, u_3\,v_3 + f\, u_3^2\, v_3  = 0 \,, \\[1.5ex]
        q_4^-&\colon (u_3,v_3)=\left(0\,,0\right) \in { \{ y_1 =0\}}  \quad \implies \quad  e = 0 \,. 
    \end{aligned}
\end{equation}
The last constraint yields 
\vspace*{-1ex}

\begin{equation}
    \begin{aligned}
        y_1&\colon  \frac{3}{2}\, f +c\,v_4 + f\, u_4^2\, v_4  = 0 \,, \\[1.5ex]
        q_5^-&\colon (u_4,v_4)=\left(0\,,\frac{3}{8}\,z\right) \in { \{ y_1 =0\}}  \quad \implies \quad  f = - \frac{c}{4}\,z \,. 
    \end{aligned}
\end{equation}
The change of variables is given by
\vspace*{-1ex}

\begin{equation}
\begin{cases}
    y_1(x_4,y_4) = B\, y_4 + C \,,\\[1.5ex] 
    x_1(x_4,y_4) = c\left( y_4^2 - \dfrac{3}{2}\, x_4 - \dfrac{z}{4} \right) \,, 
\end{cases}
\end{equation}
and to reduce the residual freedom, we impose that the variables $x_1(x_4,y_4)$, $y_1(x_4,y_4)$ are the canonical variables with respect to $H_1^{\text{P}_{\text{I\!I}}}$ 
\vspace*{-1ex}

\begin{equation}
    \begin{cases}
        x_1' = \dfrac{\partial x_1}{\partial x_4}\, x_4' + \dfrac{\partial x_1}{\partial y_4}\, y_4' + \dfrac{\partial x_1}{\partial z} = - \dfrac{\partial H_1}{\partial y_1}   \\[2ex]
        y_1' = \dfrac{\partial y_1}{\partial x_4}\, x_4' + \dfrac{\partial y_1}{\partial y_4}\, y_4' + \dfrac{\partial y_1}{\partial z} = \dfrac{\partial H_1}{\partial x_1}  
    \end{cases}
\end{equation}
where the right hand side is expressed in terms of the variables $(x_4,y_4)$. The latter are canonical with respect to $H_4^{\text{P}_{\text{I\!I}}}$ 

\vspace*{-4ex}

\begin{equation}
    \begin{aligned}
        x_4' &= - \frac{\partial H_4}{\partial y_4} \,, \qquad 
        y_1' = \frac{\partial H_4}{\partial x_4} \,. 
    \end{aligned}
\end{equation}
By equating the coefficients of similar terms between the two expressions we get 
\vspace*{-1ex}

\begin{equation}
    B = 1\,, \quad C= 0\,, \quad c = -1 \,, 
\end{equation}
giving the change of variables 
\begin{equation}
    \begin{cases}
        x_1 = \dfrac{3}{2}\,x_4 - y_4^2 + \dfrac{z}{4}\\[.7ex]
        y_1 = y_4 
    \end{cases} \,. 
\end{equation}
\vspace*{-3ex}

\section{Auxiliary function \texorpdfstring{$W$}{W}} 
\label{app:auxiliary_function}

In this appendix we exemplify how to obtain a suitable auxiliary function for a given Hamiltonian system which allows us to show that certain exceptional curves in the blown-up phase spaces cannot be reached by the flow of the vector field. Such auxiliary functions were already introduced in \cite{Shimomura2006,Shimomura2008} and \cite{halburd1, halburd2}. There, these functions were shown to satisfy a first-order differential equation with bounded coefficients, implying that these functions remain bounded themselves at any movable singularities of a solution. In this article, we employ such functions in connection with Lemma \ref{log_bounded} to show that intermediate exceptional curves cannot be reached by the solution other than at the base points, where by intermediate we mean any exceptional curve other than final exceptional curves in a cascade of blow-ups. After each blow-up, one needs to check that the auxiliary function is itself infinite on the exceptional curve, while its logarithmic derivative remains bounded there, away from the base points. Indeed, for systems with the quasi-Painlev\'e property, at any movable singularity the solution then passes horizontally through a final exceptional curve of some cascade. On the exceptional curve after the final blow-up, the auxiliary function is indeed finite. 

\subsection{\texorpdfstring{Construction of $W_3^{{\text{qsi-P}}_{\text{I\!I}}}$}{W3}}
To construct the auxiliary function, the idea is to start from the Hamiltonian, which if $z$-dependent is not itself an integral of motion, and add certain correction terms to it that ensure that this modified function remains bounded at all movable singularities. We demonstrate the process here for the Hamiltonian $H_3^{{\text{qsi-P}}_{\text{I\!I}}}$ considered in section \ref{quasi-PainleveII}. 
Noting that the total $z$-derivative of the Hamiltonian is equal to its partial $z$-derivative, we have

\vspace*{-2ex}

\begin{equation} \label{eq:H3_app_B}
\frac{d H_3^{{\text{qsi-P}}_{\text{I\!I}}}}{dz} = \frac{\partial H_3^{{\text{qsi-P}}_{\text{I\!I}}}}{\partial z} =  - c_3'(z) \,x\, y - \frac{1}{3}\, c_2'(z)\, y^3 - \frac{1}{2}\, c_1'(z)\, y^2 - c_0'(z)\, y\,.
\end{equation}
To find suitable correction terms, we consider the possible types of leading-order behaviour at movable singularities
\vspace*{-4ex}

\begin{align} 
&\begin{cases} \label{eq:behaviours_x_y_for_W_1}
y(z) = \dfrac{i}{\sqrt{2}\,(z-z_{\ast})^{1/2}} + \dfrac{i}{2\sqrt{2}}\,c_3(z_{\ast})\,(z-z_\ast)^{1/2} +\dfrac{2}{15} \, c_2(z_\ast)\,(z-z_\ast) + \mathcal{O}\big((z-z_{\ast})^{3/2}\big)\,,  \\[1.7ex]
 x(z) = -\dfrac{i}{\sqrt{2}\,(z-z_{\ast})^{3/2}}+ \dfrac{i}{2\sqrt{2}}\,\dfrac{c_3(z_\ast)}{(z-z_\ast)^{1/2}} -\dfrac{1}{15}\, c_2(z_\ast)  + \mathcal{O}\big((z-z_{\ast})^{1/2}\big)\,,
\end{cases} \\[1ex]
&\begin{cases} \label{eq:behaviours_x_y_for_W_2}
y(z) = \dfrac{1}{\sqrt{2}\, (z-z_{\ast})^{1/2}} - \dfrac{c_3(z_\ast)}{2\sqrt{2}}\,(z-z_\ast)^{1/2} + \mathcal{O}\big((z-z_{\ast})\big)\,, \\[1.7ex] 
x(z) = -\dfrac{1}{3}\,c_2(z_{\ast}) - \dfrac{c_1(z_\ast)}{\sqrt{2}}\,(z-z_\ast) + \mathcal{O}\big((z-z_{\ast})^{3/2}\big)\,. 
\end{cases}
\end{align}
Inserting these expansions for $x$ and $y$ into \eqref{eq:H3_app_B} we find an expansion of ${d H_3^{{\text{qsi-P}}_{\text{I\!I}}}}/{dz}$ in $(z-z_\ast)^{1/2}$, which after integration gives an expansion for $H_3^{{\text{qsi-P}}_{\text{I\!I}}}$ itself, which however is singular at a movable singularity $z_\ast$. By adding certain correction terms to $H_3^{{\text{qsi-P}}_{\text{I\!I}}}$ we construct a function $W$ that remains bounded at movable singularities for both types of leading order behaviours. 
By expanding each coefficient $c_k(z) = c_k(z_{\ast}) + c_k'(z_{\ast}) (z-z_{\ast}) + \cdots$, for each leading behaviour in~\eqref{eq:behaviours_x_y_for_W_1} and~\eqref{eq:behaviours_x_y_for_W_2}  we arrive at an expansion of~\!\eqref{eq:H3_app_B} in $z-z_{\ast}$. For~\eqref{eq:behaviours_x_y_for_W_1} we obtain
\begin{equation}
    \frac{d H_3^{{\text{qsi-P}}_{\text{I\!I}}}}{dz} = -\frac{1}{2} \frac{c_3'(z_{\ast})}{(z-z_{\ast})^2} +  \frac{i}{6\sqrt{2}} \frac{ c_2'(z_{\ast})}{(z-z_{\ast})^{3/2}} + \frac{1}{4} \frac{c_1'(z_{\ast})-2\,c_3''(z_{\ast})}{z-z_{\ast}} + \mathcal{O}\left((z-z_{\ast})^{-1/2}\right) \,,
\end{equation}
and for~\eqref{eq:behaviours_x_y_for_W_2} we have
\begin{equation}
     \frac{d H_3^{{\text{qsi-P}}_{\text{I\!I}}}}{dz} = -\frac{c_2'(z_\ast)}{6\sqrt{2}(z-z_\ast)^{3/2}} - \frac{c_1'(z_\ast)}{4(z-z_\ast)} + \mathcal{O} \left( (z-z_\ast)^{-1/2} \right) \,. 
\end{equation}
The numerators of the terms to power $(z-z_\ast)^{-1}$ vanish by the conditions found in the blow-up process for the system to become regular after the final blow-up. Thus,  integrating these series we find expansions for $H_3^{{\text{qsi-P}}_{\text{I\!I}}}$,
\begin{equation}
    H_3^{{\text{qsi-P}}_{\text{I\!I}}} =  \frac{1}{2(z-z_\ast)} - \frac{i}{3\sqrt{2}}\,\frac{ c_2'(z_\ast)}{(z-z_\ast)^{1/2}} + \mathcal{O}(1)
\end{equation}
and 
\begin{equation}
    H_3^{{\text{qsi-P}}_{\text{I\!I}}} = \frac{c_2'(z_\ast)}{3\sqrt{2}(z-z_\ast)^{1/2}} +  \mathcal{O}(1),
\end{equation}
respectively. To obtain the function $W$ we add terms to the Hamiltonian that eliminate all the singular (unbounded) terms, i.e.\ all terms in negative powers of $z-z_\ast$ by adding the following 
\begin{equation}
\begin{aligned}
    W_3^{{\text{qsi-P}}_{\text{I\!I}}}(z) =&~ H_3^{{\text{qsi-P}}_{\text{I\!I}}}\big(x(z),y(z);z\big) + \frac{1}{2} \frac{x(z)}{y(z)} + \frac{c_2'(z)}{6} \frac{x(z)}{y(z)^2} + \frac{c_2(z)\,c_2'(z)}{9}\,\frac{y(z)}{x(z)}  \,,
\end{aligned}
\end{equation}
obtaining the expression in \eqref{eq:W_3_QP_II}. We verified that at each step of the blow-up procedure the function $W_3^{{\text{qsi-P}}_{\text{I\!I}}}$ is bounded on the exceptional divisors except for the base points, and it is finite in correspondence of the last exceptional curve for both the branches  exhibiting the behaviours \eqref{eq:behaviours_x_y_for_W_1} and \eqref{eq:behaviours_x_y_for_W_2} respectively.  

\subsection{\texorpdfstring{Construction of $W_1^{{\text{qsi-P}}_{\text{I\!V}}}$}{W1}}
We proceed with $W_1^{{\text{qsi-P}}_{\text{I\!V}}}$ similar to the case above, only that there are three types of leading-order behaviour~\eqref{leading_QPIV_H1_y} to accommodate, of the form 
\begin{equation} \label{eq:leading_order_W1P4_app}
    \begin{cases}
        x(z) \sim \dfrac{2}{z-z_\ast},\\[2ex] 
        y(z) \sim a_0 (z-z_\ast),
    \end{cases} \hspace{5ex} \begin{cases}
        x(z) \sim \dfrac{1}{2(z-z_\ast)},\\[2ex] 
        y(z) \sim \dfrac{\sigma}{(z-z_\ast)^{1/2}}, \text{ with } \sigma \in \{1,i\}\,.
    \end{cases}
\end{equation}
We note that for the singularities of the first leading order in \eqref{eq:leading_order_W1P4_app} no additional terms are required to render $W_1^{{\text{qsi-P}}_{\text{I\!V}}}$ bounded. We therefore focus on the two remaining leading orders in \eqref{eq:leading_order_W1P4_app}. The terms we add to the Hamiltonian to obtain a bounded function are of the form $b_k(z) \,x(z)\,y(z)^{-k}$, for $k=0,1,2,\dots$, namely
\begin{equation}
\begin{aligned}
   W_1^{{\text{qsi-P}}_{\text{I\!V}}} &= H_1^{{\text{qsi-P}}_{\text{I\!V}}} + b_0(z) \, x(z) + b_1(z) \frac{x(z)}{y(z)} \,, \\[1ex]
   b_0(z)&=\frac{1}{2}\,, \qquad b_1(z) = \frac{a_3'(z)}{3} - \frac{3z}{2}\,. 
\end{aligned}
\end{equation}
We note however, that these terms now conflict with the bounded behaviour of $W$ for the third leading-order. To remedy this, we note that we can insert additional bounded terms $b_k(z) \,x(z)\,y(z)^{-k}$, $k \geq 2$ without changing the behaviour of $W$. We add infinitely many such terms to complete this sum into a geometric series in $u=y^{-1}$, the sum of which is a bounded quantity also for the third possible leading-order. Since the two coefficients, $b_0$ and $b_1$, are fixed, the expression for the sum of this series can be taken to be of the form 
\begin{equation}
\begin{aligned} 
    \frac{A_0 \,x(z)}{1 + A_1(z)\,u} &= A_0(z)\,x(z) - A_0\, A_1(z)\, x(z)\, u + \cdots, \\[1ex]
    A_0 &= b_0 = \frac{1}{2}\,, \qquad A_1(z) = 3z - \frac{2}{3} \,a_3'(z)\,.
\end{aligned}
\end{equation}
One can now check, using computer algebra, that in each intermediate step of each of the three cascades of blow-ups, the function
\begin{equation}
    W_1^{{\text{qsi-P}}_{\text{I\!V}}} = H_1^{{\text{qsi-P}}_{\text{I\!V}}} + \frac{1}{2} \frac{x(z)\, y(z)}{y(z) + A_1(z)}
\end{equation}
is infinite on the exceptional curve while the logarithmic derivative of $W_1^{{\text{qsi-P}}_{\text{I\!V}}}$ is bounded. Furthermore, on the final exceptional curve of each cascade, $W_1^{{\text{qsi-P}}_{\text{I\!V}}}$ is bounded.

The corresponding bounded function for the Hamiltonian system derived from $H_2^{{\text{qsi-P}}_{\text{I\!V}}}$ is very similar, 
\begin{equation}
    W_2^{{\text{qsi-P}}_{\text{I\!V}}} = H_2^{{\text{qsi-P}}_{\text{I\!V}}} + \frac{1}{2} \frac{x(z)\, y(z)}{y(z) + B_1(z)}, \qquad B_1(z) = z - \frac{2}{3} \,a_3'(z)\,.
\end{equation}

\bibliographystyle{style}
\bibliography{thebiblio}

\end{document}